\documentclass[11pt,a4paper,oneside]{article}
\usepackage{amsmath,amsthm}
\usepackage{amsfonts}
\usepackage{graphicx, color}
\usepackage{amssymb,dsfont,hyperref,color}
\usepackage[left=2.5cm,right=2.5cm,top=3cm,bottom=3cm]{geometry}
\allowdisplaybreaks

\newtheorem{theorem}{Theorem}
\newtheorem{proposition}[theorem]{Proposition}
\newtheorem{corollary}{Corollary}

\newtheorem{lemma}{Lemma}
\newtheorem{definition}{Definition}
\newtheorem{remark}{Remark}

\numberwithin{equation}{section}

\newcommand{\vertiii}[1]{{\left\vert\kern-0.25ex\left\vert\kern-0.25ex\left\vert #1
    \right\vert\kern-0.25ex\right\vert\kern-0.25ex\right\vert}}

\DeclareMathOperator{\argmin}{argmin\,}

\DeclareMathOperator{\pen}{pen}

\newcommand{\hmin}{h_{\min}}
\newcommand{\hatfhat}{\hat f_{\hat h }}

\def\S{\mathbb{S}^{d-1}}

\def\P{\mathbb{P}}
\def\R{\mathbb{R}}
\def\E{\mathbb{E}}

\def\L{\mathbb{L}}
\def\Kh{K_{h^2}}



\title{Adaptive optimal kernel density estimation for directional data}
\author{Thanh Mai Pham Ngoc
 \thanks{Laboratoire de Math\'ematique d'Orsay, Univ. Paris-Sud, CNRS,  Universit\'e Paris-Saclay, 91405 Orsay, France.  Email: thanh.pham\_ngoc@math.u-psud.fr}}

\date{\today}
\begin{document}
\maketitle

\begin{abstract}
We focus on the nonparametric density estimation problem with directional data. We propose a new rule for bandwidth selection for kernel density estimation. Our procedure is automatic, fully data-driven and adaptive to the smoothness degree of the density. We obtain an oracle inequality and optimal rates of convergence for the $\mathbb{L}_2$ error. Our theoretical results are illustrated with simulations.
\end{abstract}



\section{Introduction }

Directional data arise in many fields such as wind direction for the circular case, astrophysics, paleomagnetism, geology for the spherical case.  Many efforts have been put to devise statistical methods to tackle the density estimation problem. We refer to \cite{Mardia} and more recently to \cite{Ley} for a comprehensive view. Nonparametric procedures have been well developed. In this article we focus on kernel estimation but we may cite a series of works using projection methods on localized bases adapted to the sphere (\cite{Baldi}, \cite{Tani}).  Classical references for kernel estimation with directional data include the seminal papers of \cite{Hall} and \cite{Bai}. It is well-known that the choice of the bandwidth is a key and intricate issue when using kernel methods. Various techniques for selecting the bandwidth have been suggested since the popular cross-validation rule by \cite{Hall}. We shall cite plug-in and refined cross-validatory methods in \cite{Taylor} and \cite{Oliveira} for the circular case, \cite{DiMarzio} on the torus. More recently, \cite{Garcia2} devised an equivalent of the rule-of-thumb of \cite{Silverman} for directional data, whereas \cite{Amiri} explored computational problems with recursive kernel estimators based on the cross-validation procedure of \cite{Hall}. But to the best of our knowledge, all the rules proposed so far for selecting the bandwidth are empirical. Although they prove efficient in practice, only little attention has been put on their theoretical properties. \cite{Klemela2} studied convergence rates for $\mathbb{L}_2$ error over some regularity classes for the kernel estimator for the estimation of the density and its derivatives but the asymptotically optimal bandwidth depends on the density and its smoothness degree which is unfeasible for applications. In the present paper, we try to fill the gap between theory and practice. Our goal is two-fold. We aim at devising an automatic and fully data driven choice of the kernel bandwidth so that the resulted estimator achieves optimal rates in the minimax sense for the $\mathbb{L}_2$ risk over some regularity classes. We emphasize that the estimator is adaptive to the smoothness degree of the underlying density. It means that the method does not require the specification of the regularity of the density. Our work is inspired by very recent techniques developped for the multivariate case.   In the problem of multivariate kernel density estimation, adaptive minimax approaches have been tackled in  a remarkable series of papers by \cite{Lepski1}, \cite{Lepski2}, \cite{Gold} and very recently by \cite{Lacour}.  Our methodology is inspired from the PCO (Penalized Comparison of Overfiting) procedure of \cite{Lacour}. It is based on concentration inequalities for $U-$statistics. Last but not least, our procedure is simple to be implemented and in examples based on simulations, it shows quite good performances in a reasonable computation time.

This paper is organized as follows. In section \ref{estimation}, we present our estimation procedure. In section \ref{vitesses} we provide an oracle inequality and rates of convergences of our estimator for the MISE. Section \ref{simus} gives some numerical
illustrations. Section \ref{preuves} gives the proofs of theorems. Finally the Appendix gathers some technical lemmas.

\vspace{0.5cm}

\textbf{Notations} For two integers $a$, $b$, we denote $a\wedge b := \min(a, b)$ and $a\vee b := \max(a, b)$. And  $\lfloor y \rfloor $denotes the largest integer smaller than $y$ such that $  \lfloor y \rfloor \leq  y < \lfloor y \rfloor+1$.

 Depending on the context, $ \| \cdot \|$ denotes the classical $\mathbb{L}_2$-norm on $\R$ or $\S$, $\langle \cdot, \cdot \rangle$ the associated scalar product. For a vector $x\in \R^d$, $\| x \|$ stands for the Euclidian norm on $\R^d$. And $ \| \cdot \|_{\infty} $ is the usual  $\mathbb{L}_{\infty}$-norm on $\S$. 
  
The scalar product of two vectors $x$ and $y$, is denoted by $x^T y$, where $^T$ is the transpose operator. 

\section{Estimation procedure}\label{estimation}

We observe $n$ i.i.d observations $X_1, \dots, X_n$ on $\S$ the unit sphere of $\R^{d}$, $d \geq 3$. The $X_i$'s are absolutely continuous with respect to the Lebesgue measure $\omega_d$ on $\S$ with common density $f$. Therefore, a directional density $f$ satisfies
$$\int_{\S} f(x) \omega_d(dx)=1. $$
 We aim at constructing an adaptive kernel estimator of the density $f$ with a fully data-driven choice of the bandwidth. 

\subsection{Directional approximation kernel}

We present some technical conditions that are required for the kernel.

\textbf{Assumption 1} The kernel $K: [0, \infty) \rightarrow [0, \infty) $ is a bounded  and Riemann integrable function such that 
\[
0 < \int_{0}^{\infty} x^{(d-3)/2}K(x)dx < \infty. 
\]

Assumption 1 is classical in kernel density estimation with directional data, see for instance Assumptions D1-D3 in \cite{Garcia} and Assumption A1 in \cite{Amiri}. An example of kernel which satisfies Assumption 1 is the popular von-Mises kernel $K(x)=e^{-x}$.

\vspace{0.5cm}
Now following \cite{Klemela2} we shall define what is called a  kernel of class $s$. Let
\[
\alpha_i(K)= \int_0^{\infty} x^{(i+d-3)/2} K(x) dx, \quad i \in \mathbb{N}.
\]

\textbf{Assumption 2}
Let $s \geq 0$ be even. The kernel $K$ is of class $s$ i.e it is a measurable function $K: [0, +\infty[ \rightarrow \R$ which satisfies :
\begin{itemize}
\item $\alpha_i (K) < +\infty$ for $i=0, \dots, s$
\item $\alpha_0(K) \neq 0$,
\item $ \int_{0}^{h^{-2}} x^{(2i+d-3)/2} K(x) dx = o(h^{s-2i})$ for $i= 1, \dots, s/2-1$, when $h$ tends to $0$.
\end{itemize}

\subsection{Family of directional kernel estimators}

We consider the following classical directional kernel density estimator  $\left( K_h(x,y) := K \left (\frac{1-x^Ty}{h}  \right )\right )$
$$\hat f_h(x)= \frac{c_0(h)}{n} \sum_{i=1}^n K\left (\frac{1-x^TX_i}{h^2} \right )=  \frac{c_0(h)}{n} \sum_{i=1}^n K_{h^2}\left ({x,X_i} \right ), \quad x \in \S,$$
where $K$ is a kernel satisfying Assumption 1 and  $c_0(h)$ a normalizing constant such that $\hat f_h(x)$ integrates to unity:
\begin{equation}\label{c_0}
 c_0^{-1}(h):= \int_{\S} K_{h^2}(x,y)\omega_d(dy).
\end{equation}
It remains to select a convenient value for $h$.

\subsection{Bandwidth selection}
In kernel density estimation, a delicate step consists in  selecting the proper bandwith $h$ for $\hat f_h$.  We suggest the following data-driven choice of bandwith $\hat h$ inspired from \cite{Lacour}. We name our procedure SPCO (Spherical Penalized Comparison to Overfitting). We denote $  \E(\hat f_h):=f_h.$ Our selection rule is the following: 

\begin{equation}\label{hhat}
\hat h := \argmin_{h \in \mathcal{H}} \left \{ \| \hat f_h -\hat f_{h_{\min}} \|^2 +\textrm{pen}_\lambda(h) \right \} \quad \lambda \in  \R,
\end{equation}
where $h_{\min} = \min \mathcal{H}$ 
\begin{equation}\label{penalite}
\textrm{pen}_\lambda(h)  := \frac{\lambda c_0^2(h)c_2(h)}{n} -\frac 1 n \int_{\S} \left( c_0(h_{\min}) K_{h_{\min}^2}(x,y) - c_0(h) K_{h^2}(x,y)  \right )^2 \omega_d(dy),
\end{equation}


\begin{equation*}
c_2(h):= \int_{\S} K^2_{h^2}(x,y)\omega_d(dy),
\end{equation*}
and $\mathcal H$ a set of bandwiths  defined by 
\begin{equation}\label{H}
\mathcal{H}=\left  \{h:  \left ( \frac{\| K\|_{\infty}}{n}\frac{1}{R_0(K)}\right )^{\frac{1}{d-1}} \leq h \leq 1,\textrm{ and } h^{-1} \;\textrm{is an integer} \right \},
\end{equation}
with  $ R_0(K) := 2^{(d-3)/2}\sigma_{d-2}\alpha_0(K)$ and $\sigma_{d-1}= \omega_{d}(\S)$ denotes the area of $\S$. We recall that
$
\sigma_{d-1} = \frac{2\pi^{d/2}}{\Gamma(d/2)},
$
with $\Gamma$ the Gamma function. 

\bigskip
The estimator of $f$ is $\hat f_{\hat h}$.

\bigskip

The procedure SPCO involves a real parameter $\lambda$. In section  \ref{vitesses} we study how to choose the optimal value of $\lambda$ leading to a fully data driven procedure.

\begin{remark}\label{remarque1}

Note that $c_0(h)$, $c_2(h)$ and $\textrm{pen}_{\lambda}(h)$ do not depend on $x$.  Indeed its is known (see $\cite{Hall}$) that if $y$ is a vector and $x$ a fixed element of $\S$, then denoting $t= x^Ty$ their scalar product, we may always write 
$$
y= t x +(1-t^2)^{1/2} \xi,
$$

where $\xi$ is a unit vector orthogonal to $x$. Further, the area element on $\S$ can be written as
$$
\omega_d(dx)=(1-t^2)^{\frac{d-3}{2}} dt \; \omega_{d-1}(d\xi ).
$$
Thus, using these conventions, one obtains
\begin{eqnarray*}
c_0^{-1}(h)&=&  \int_{\S} K \left (\frac{1-x^Ty}{h^2} \right )\omega_d(dy)  \\
&=& \int_{\mathbb{S}^{d-2}} \int_{-1}^1 K \left(\frac{1- x^T( t x +(1-t^2)^{1/2} \xi)}{h^2}\right)(1-t^2)^{\frac{d-3}{2}} dt \; \omega_{d-1}(d\xi) \\
&= &   \int_{\mathbb{S}^{d-2}} \omega_{d-1}(d\xi)  \int_{-1}^1 K \left(\frac{1- t x^Tx}{h^2} \right) (1-t^2)^{\frac{d-3}{2}} dt \\
&= & \sigma_{d-2}\int_{-1}^1 K \left(\frac{1- t }{h^2} \right) (1-t^2)^{\frac{d-3}{2}} dt.
\end{eqnarray*}
Using similar computations,  one gets that 
$$
c_2(h) = \sigma_{d-2}\int_{-1}^{1} K^2\left (\frac{1-t}{h^2} \right )(1-t^2)^{(d-3)/2} dt,
$$
and 
\[
\textrm{pen}_\lambda(h)= \frac{\lambda c_0^2(h)c_2(h)}{n} - \frac{ \sigma_{d-2}}{n} \int_{-1}^1 \left( c_0(h_{\min}) K \left (\frac{1-t}{h_{\min}^2} \right ) -  c_0(h) K \left (\frac{1-t}{h^2} \right )  \right )^2 (1-t^2)^{(d-3)/2} dt.
\]
\end{remark}


\section{Rates of convergence}\label{vitesses}

\subsection{Oracle inequality}\label{section-oracle}

First, we state an oracle type inequality which highlights the bias-variance decomposition
of the $\mathbb{L}_2$ risk. $| \mathcal{H}|$ denotes the cardinality of $\mathcal{H}$. We recall that we denote $f_h := \E(\hat f _h ) $.

\begin{theorem}\label{theorem 2}
Assume that kernel $K$ satisfies Assumption 1 and $\|f \|_{\infty} < \infty$. Let  $x \geq 1$ and $\varepsilon \in (0,1)$. 
Then there exists $n_0$ independent of $f$, such that  for $n\geq n_0$  with probability larger than $1- C_1 |\mathcal H | e^{-x}$,
\begin{equation}\label{oracle}
\begin{split}
\| \hat f_{\hat h}- f\|^2 & \leq C_0(\varepsilon, \lambda) \min_{h \in \mathcal{H}} \| \hat f_h -f \|^2  \\
 & \quad + C_2(\varepsilon, \lambda)  \| f_{\hmin} -f \|^2 +C_3(\varepsilon, K,\lambda) \left ( \frac{\| f\|_{\infty}x^2}{n} + \frac{ c_0(\hmin)x^3 } {n^2 }\right ),
\end{split}
\end{equation}
where $C_1$ is an absolute constant and $C_0(\varepsilon, \lambda) =\lambda +\varepsilon$ if $\lambda \geq 1$, $C_0(\varepsilon)=\frac 1 \lambda +\varepsilon$ if $ 0 < \lambda < 1$. The constant $C_2(\varepsilon, \lambda)$ only depends on $\varepsilon$ and $\lambda$ and $C_3(\varepsilon, K,\lambda)$ only depends on $\varepsilon$, $K$ and $\lambda$.
\end{theorem}

This oracle inequality bounds the quadratic risk of SPCO
estimator by the infimum over $\mathcal{H}$ of the tradeoff between the approximation term $  \| f_{\hmin} -f \|^2$ and the variance term $ \| \hat f_h -f \|^2$. The terms $C_3(\varepsilon, K,\lambda)\left  ( \frac{\| f\|_{\infty}x^2}{n} + \frac{ c_0(\hmin)x^3 } {n^2 } \right ) $ are remaining terms. Hence, this oracle inequality justifies our selection rule. For further details about oracle inequalities and model selection see \cite{Massart}.
\bigskip 

The next theorem shows that  we cannot choose $\lambda$ too small ($\lambda<0$) at the risk of selecting a bandwidth close to $\hmin$ with high probability.  This would lead to an overfitting estimator. To this purpose, we suppose
\begin{equation}\label{biais_negligeable}
\| f-f_{\hmin}\|^2 \frac{n}{c_0^2(\hmin) c_2(\hmin)}=o(1).
\end{equation} 
This assumption is quite mild.  Indeed the variance of $\hat f_h$ is of order $\frac{c_0^2(h)c_2(h)}{n}$, thus this assumption means that the smallest bias is negligible with respect to the corresponding integrated variance. Last but not least, because $c_0^2(\hmin) c_2(\hmin)$ is of order $\hmin^{1-d}$ (see Lemme \ref{proprietes}), this assumption amounts to
$$
\| f-f_{\hmin}\|^2 = o\left( \frac{\hmin^{1-d}}{n}\right).
$$

\begin{theorem}\label{theorem_petit}
Assume that kernel $K$ satisfies Assumption 1 and $\|f \|_{\infty} < \infty$. Assume also (\ref{biais_negligeable}) and
\[
\frac{\| K\|_{\infty}}{n}\frac{1}{R_0(K)} \leq h_{\min}^{d-1} \leq \frac{(\log n)^{\beta}}{n}, \quad \beta >0. 
\]
Then if we consider $\textrm{pen}_\lambda(h)$ defined in (\ref{penalite}) with $\lambda <0$, then
we have for $n$ large enough, with probability larger than $1-C_1|\mathcal{H}| e^{-(n/log n)^{1/3}}$:
\[
\hat h \leq C(\lambda) \hmin \leq  C(\lambda)  \left (\frac{(\log{n})^{\beta}}{n} \right)^{\frac{1}{d-1}},
\]
where $C_1$ is an absolute constant and $C(\lambda)= \left( 1.23 \left( 2.1-\frac 1 \lambda \right ) \right)^{\frac{1}{d-1}}$.

\end{theorem}

\begin{remark}\label{remark_lambda}
Theorem \ref{theorem_petit} invites us to discard $\lambda <0$. Now considering oracle inequality (\ref{oracle}), $\lambda=1$ yields the minimal value of the leading constant $C_0(\varepsilon, \lambda)=\lambda + \varepsilon $ . Thus, the theory urges us to take the optimal value $\lambda=1$ in the SPCO procedure.  Actually, we will see in the numerical section that the choice $\lambda=1$ is quite efficient. 
\end{remark}

\subsection{Tuning-free estimator and rates of convergence }\label{final-estimateur}

Results of Section \ref{section-oracle} about the optimality of $\lambda=1$ enable us to devise our  tuning-free estimator $\hat f_{\check h}$ with bandwidth $\check h$ defined as follows

\begin{equation}\label{hhat-tuning-free}
\check h := \argmin_{h \in \mathcal{H}} \left \{ \| \hat f_h -\hat f_{h_{\min}} \|^2 +\textrm{pen}(h) \right \}
\end{equation}
and
\begin{equation}\label{penalite-free}
\textrm{pen}(h)  := \frac{ c_0^2(h)c_2(h)}{n} -\frac 1 n \int_{\S} \left( c_0(h_{\min}) K_{h_{\min}^2}(x,y) - c_0(h) K_{h^2}(x,y)  \right )^2 \omega_d(dy),
\end{equation}

and  $\hmin = \min \mathcal H  $  where $ \mathcal H$ is a set of bandwiths  defined by 
\begin{equation}\label{H}
\mathcal{H}=\left  \{h:  \left ( \frac{\| K\|_{\infty}}{n}\frac{1}{R_0(K)}\right )^{\frac{1}{d-1}} \leq h \leq 1,\textrm{ and } h^{-1} \;\textrm{is an integer} \right \}.
\end{equation}

The following corollary of Theorem \ref{theorem 2} states an oracle inequality satisfied by $\hat f_{\check h }$. It will be central to compute rates of convergence. 

\begin{corollary}\label{corollaire}
Assume that kernel $K$ satisfies Assumption 1 and $\|f \|_{\infty} < \infty$. Let  $x \geq 1$ and $\varepsilon \in (0,1)$. 
Then there exists $n_0$ independent of $f$, such that  for $n\geq n_0$  with probability larger than $1- C_1 |\mathcal H | e^{-x}$,
\begin{equation}\label{oracle}
\begin{split}
\| \hat f_{\check h}- f\|^2 & \leq C_0(\varepsilon) \min_{h \in \mathcal{H}} \| \hat f_h -f \|^2  \\
 & \quad + C_2(\varepsilon)  \| f_{\hmin} -f \|^2 +C_3(\varepsilon, K) \left ( \frac{\| f\|_{\infty}x^2}{n} + \frac{ c_0(\hmin)x^3 } {n^2 }\right ),
\end{split}
\end{equation}
where $C_1$ is an absolute constant and $C_0(\varepsilon) =1 +\varepsilon$. The constant $C_2(\varepsilon)$ only depends on $\varepsilon$  and $C_3(\varepsilon, K)$ only depends on $\varepsilon$ and $K$.
\end{corollary}

We now compute rates of convergence for the MISE (Mean Integrated Square Error) of our estimator $\hat f_{\check h}$ over some smoothness classes.
\cite{Klemela2} defined suitable smoothness classes for the study of the MISE. In particular, theses regularity classes involve a concept of an "average" of directional derivatives which was first defined in \cite{Hall}. Let us recall the definition of these smoothness classes. 
\bigskip

  Let $ \eta \in \S$, $d \geq 3$ and $T_\eta  = \{ \xi \in \S | \xi \perp \eta \}$. Let $\phi_\eta : \S \setminus \{ \eta, -\eta\} \rightarrow T_\eta \times ]0, \pi [  $ be a parameterization of $\S$ defined by
$$
\phi_\eta^{-1}(\xi, \theta)= \eta \cos(\theta) +\xi \sin(\theta).
$$

When $g: \R^{d} \rightarrow \R$ and $x, \xi \in \R^{d}$, define the derivative of $g$  at $x$ in the direction of $\xi$ to be $D_{\xi} g(x)= \lim_{h \rightarrow 0} h^{-1}[g(x+h\xi)-g(x)]$ and $D^{s}_{\xi} g = D_{\xi} D_{\xi}^{s-1} g,$ for $s \geq 2$ an integer. 

We now shall define the derivative  of order $s$.
\begin{definition}
Let $f : \S \rightarrow \R $ define $\bar f : \R^{d} \rightarrow \R$ by $\bar f(x)= f(x/ \| x\|)$. The derivative of order $s$ is  $D^s f : \S \rightarrow \R$ defined by
$$
D^s f(x) = \frac{1}{\sigma_{d-1}} \int_{T_x} D^s_{\xi} \bar f(x) \omega_d(d\xi),
$$
 where $d \geq 3$, $T_x= \{ \xi \in \S : \xi \perp x \}$.
 \end{definition}


\begin{definition}
When $f: \S \rightarrow \R$ define $\tilde D^s f : \S \times \R \rightarrow \R$ by
$$
\tilde D^s f(x, \theta) = \frac{1}{\sigma_{d-1}} \int_{T_x} D^s_{\phi_x^{-1}(\xi, \theta + \frac {\pi}{2})} \bar f(\phi^{-1}_x(\xi, \theta)) \omega_d(d\xi).
$$

\end{definition}

We are now able to define the smoothness class $\textbf{F}_2(s)$ (see \cite{Klemela2}). 

\begin{definition}\label{Sobolev}
Let $s \geq 2$ be even and $1 \leq p \leq \infty $. Let $\textbf{F}_2(s)$ be the set of such functions $f: \S \rightarrow \mathbb{R}$ that 
$\| D^i f \| < \infty$ for $i=0, \dots, s$ for all $x \in \S$ and for all $\xi \in T_x$, $\frac{\partial^s}{\partial \theta^s}f(\phi_{x}^{-1}(\xi, \theta))$ is continuous as a function of $\theta \in \R$, $\| \tilde D^s f(\cdot, \theta)\|$ is bounded for $\theta \in [0, \pi]$ and $\lim_{\theta \rightarrow 0} \| \tilde D^s f(\cdot, \theta) -D^s f\|= 0$.
\end{definition}

An application of the oracle inequality in Corollary \ref{corollaire} allows us to derive rates of convergence for the MISE of $\hat f_{\check h}$.



\begin{theorem}\label{vitesse}
Consider a kernel $K$ satisfying Assumption 1 and Assumption 2. For $B>0$, let us denote $\tilde{\textbf{F}}_2(s,B)$ the set of densities bounded by $B$ and belonging to $\textbf{F}_2(s)$.  Then we have 
\[
\sup_{f \in \tilde{\textbf{F}}(s,B)} \E[ \| \hat f_{\check h} -f\|^2]  \leq C(s, K, d, B) n^{\frac{-2s}{2s+(d-1)}},
\]
with $C(s,K,d,B)$ a constant depending on $s$, $K$, $d$ and $B$.
\end{theorem}

 Theorem \ref{vitesse} shows that the estimator $ \hat f_{\check h}$ achieves the optimal rate of convergence in dimension $d-1$ (minimax rates for multivariate density estimation are studied in \cite{Ibra1}, \cite{Ibra2}, \cite{Ibra3}). Furthermore, our statistical procedure is adaptive to the smoothness $s$. It means that it does not require the specification of $s$. 

\section{Numerical results}\label{simus}
We  investigate the numerical performances of our fully data-driven estimator $\hat f_{\check h }$ defined in section \ref{final-estimateur} with simulations. We compare  $\hat f_{\check h}$   to the widely used cross-validation estimator and to the "oracle" (that will be defined later on). We focus on the unit sphere $\mathbb{S}^2$ (i.e the case $d=3$).

\begin{figure}[!h]\begin{center}
\begin{tabular}{c}
\includegraphics[scale=0.35]{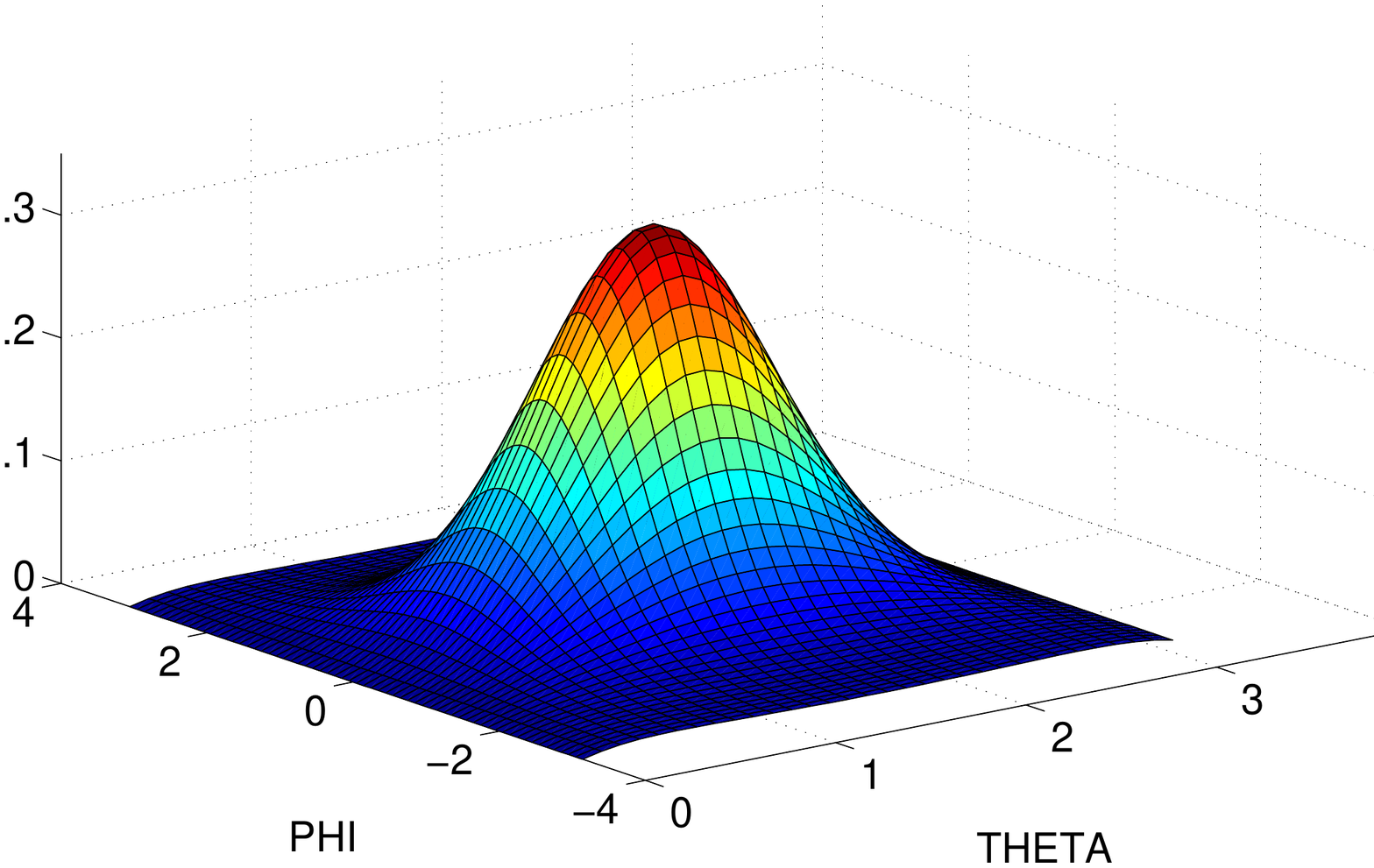} 
\end{tabular}
\caption{The density $f_{1,vM}$ in spherical coordinates}
 \label{densite-fisher}
\end{center}
\end{figure}

We aim at estimating the  von-Mises Fisher density  $f_{1,vM}$ (see Figure \ref{densite-fisher}): 
$$
f_{1,vM}= \frac{\kappa}{2\pi (e^{\kappa}-e^{-\kappa})}e^{\kappa x^T\mu},
$$
with $\kappa=2$ and  $\mu= (1,0,0)^T$. We recall that $\kappa $ is the concentration parameter and $\mu$ the directional mean. 
We also estimate the mixture of two von-Mises Fisher densities, $f_{2,vM}$:
$$
f_{2,vM}=\frac 4 5 \times   \frac{\kappa}{2\pi (e^{\kappa}-e^{-\kappa})} e^{\kappa x^T\mu} + \frac 1 5  \times \frac{\kappa'}{2\pi (e^{\kappa'}-e^{-\kappa'})} e^{\kappa' x^T\mu'},
$$
with $\kappa'=0.7$ and $\mu'=(-1,0,0)^T$. Note that the smaller the concentration parameter is, the closer to the uniform density the  von-Mises Fisher density is. 

 Now let us define what the "oracle" $ \hat f_{ h_{oracle}}$ is. The bandwidth $h_{oracle}$ is defined as 
\[
h_{oracle}=\argmin_{h\in \mathcal{H}}   \| \hat f_h -f\|^2. 
\]
$h_{oracle}$ can be viewed as  the "ideal" bandwidth since it uses the specification of the density of interest $f$ which is here either $f_{1,vM}$ or $f_{2, vM}$.  Hence, the performances of $\hat f_{ h_{oracle}}$ are used as a benchmark.

In the sequel we present detailed results for $f_{1,vM}$, namely risk curves and graphic reconstructions and we compute MISE both for $f_{1,vM}$ and $f_{2,vM}$. We use the von-Mises kernel $K(x)=e^{-x}$.
\medskip

\begin{figure}[!h]\begin{center}
\begin{tabular}{cc}
\includegraphics[scale=0.35]{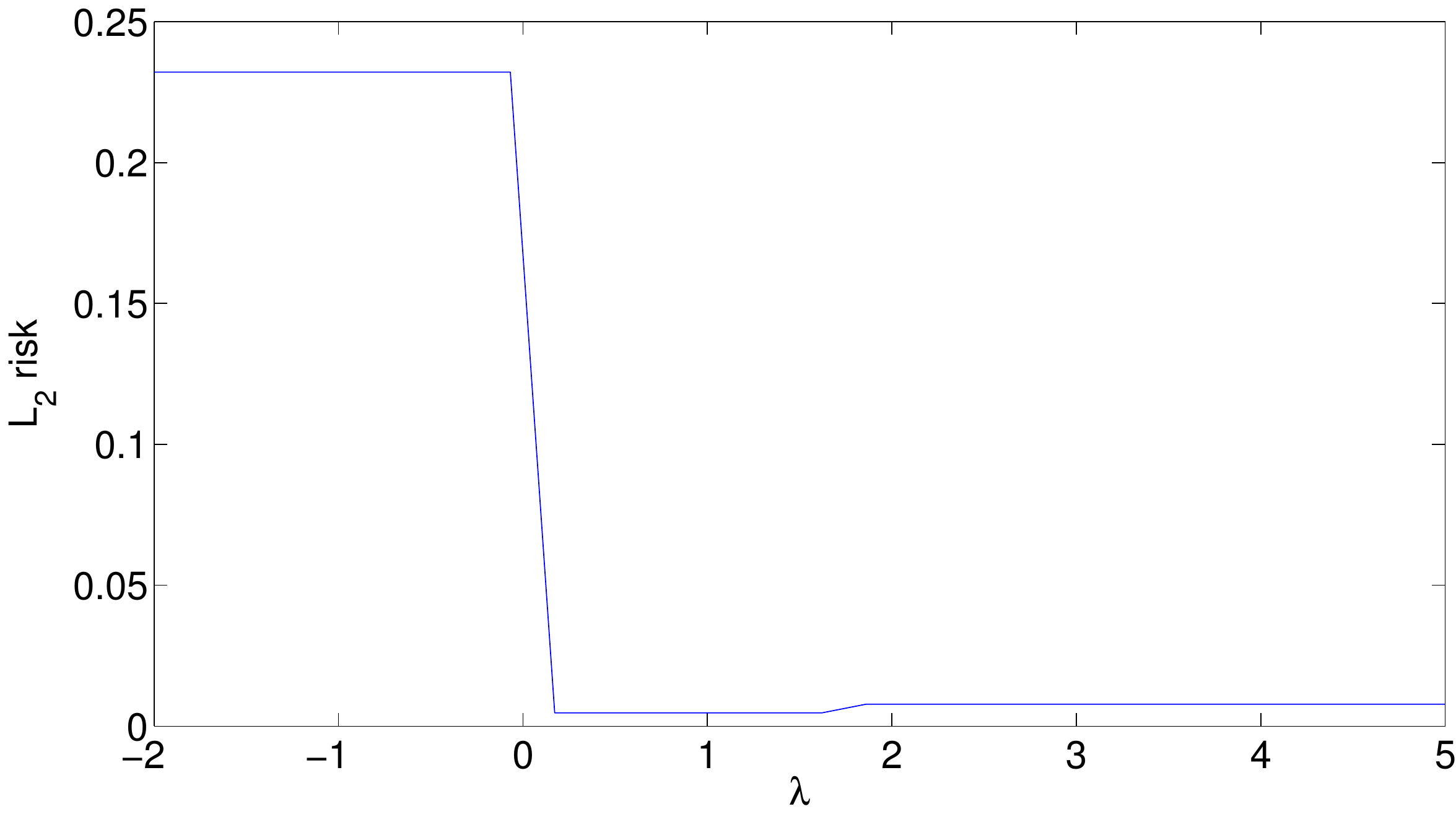} &
\includegraphics[scale=0.35]{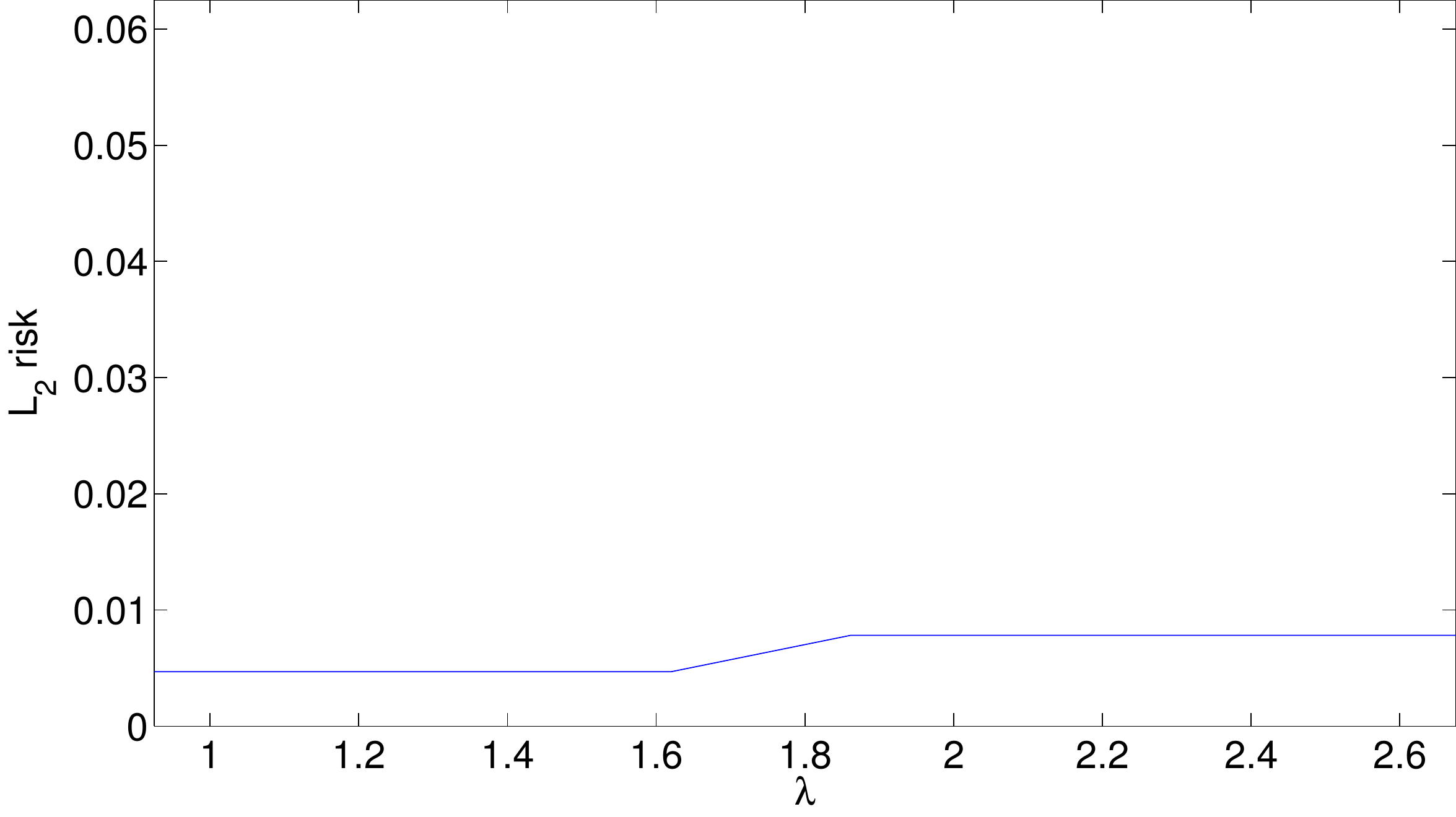} \\
a & b 
\end{tabular}
\caption{a/ $\mathbb{L}_2$ risk of $\hat f_{\check h}$ to estimate $f_{1,vM}$ in function of parameter $\lambda$  \quad  b/ A zoom } 
 \label{sautdim}
\end{center}
\end{figure}

Before presenting the performances of the various procedures, we shall remind that theoretical results of section \ref{section-oracle} have shown that setting $\lambda=1$ in the SPCO algorithm was optimal. We would like to show how simulations actually support this conclusion. Indeed, Figure \ref{sautdim} displays the $\mathbb{L}_2$-risk of $\hat f_{\hat h}$ in function of parameter $\lambda$.  Figure \ref{sautdim} a/ shows a  "dimension jump" and that the minimal risk is reached in a stable zone around $\lambda=1$: negative values of $\lambda$ lead to an overfitting estimator ($\hmin$ is chosen) with an explosion of the risk, whereas large values of $\lambda$ make the risk increase again (see a zoom on Figure \ref{sautdim} b/). Next, we will realize that $\lambda=1$ yields quite good results.
\medskip

In Lemma \ref{SPCO-pratique} of the Appendix, we clarify  the expression (\ref{hhat-tuning-free}) to be minimized to implement our estimator $\hat f_{\check h }$. We now recall the cross-validation criterion of \cite{Hall}. Let 
\[
\hat f_{h,i}(x)= \frac{c_0(h)}{n-1} \sum_{j \neq i}^n e^{-(1-^txX_j)/h^2} 
\]
then
\[
CV_2(h)= \| \hat f_h\|^2 -\frac 2 n \sum_{i=1}^n \hat f_{h,i}(x).
\]
$CV_2$ is an unbiased estimate of the MISE of $\hat f_h$. The cross-validation procedure to select the bandwidth $h$ consists in minimizing $CV_2$ with respect to $h$. We call this selected value $h_{CV_2}$.

\begin{figure}[!h]\begin{center}
\begin{tabular}{ccc}
\includegraphics[scale=0.22]{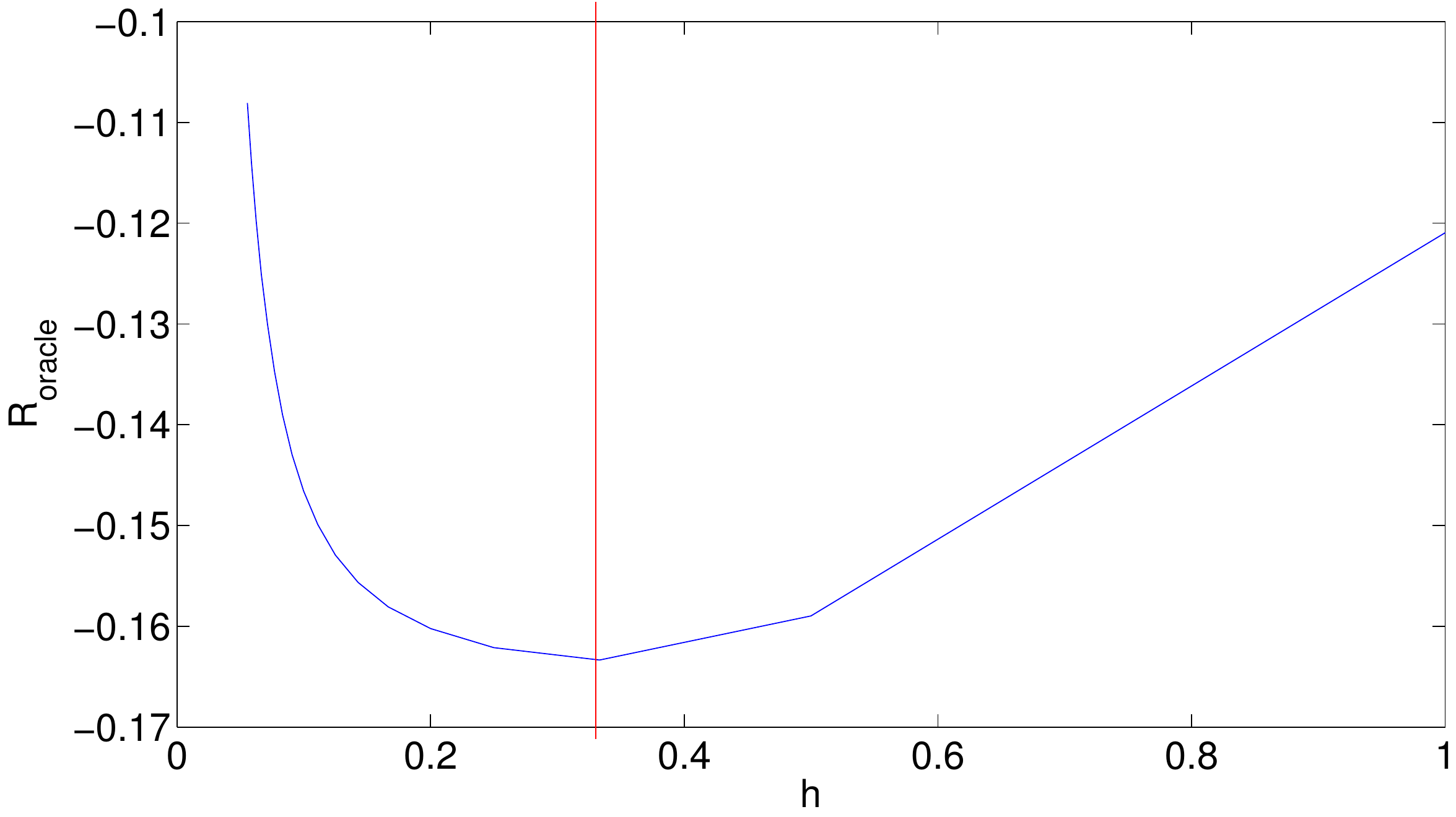}&
\includegraphics[scale=0.22]{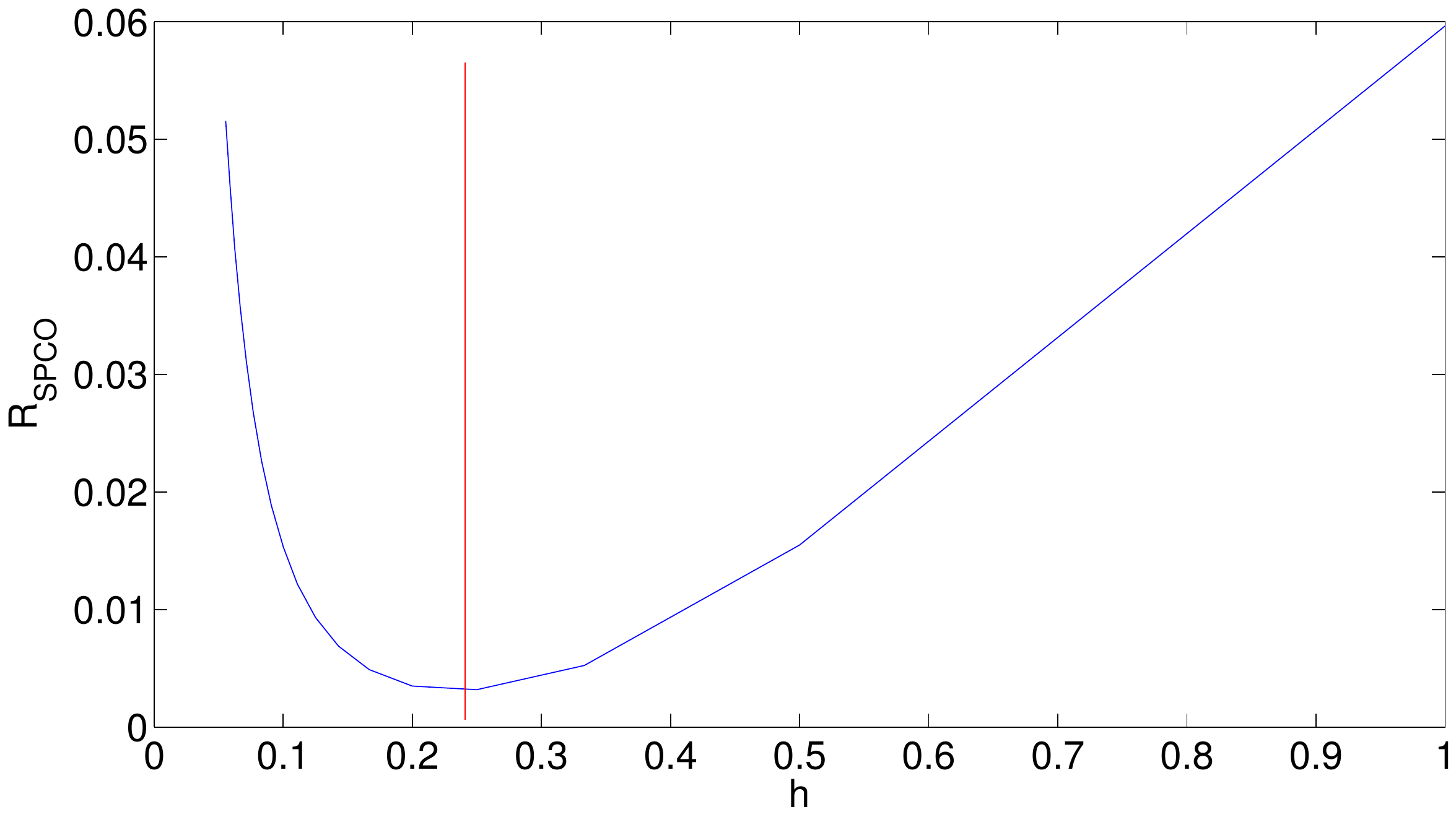}&
\includegraphics[scale=0.22]{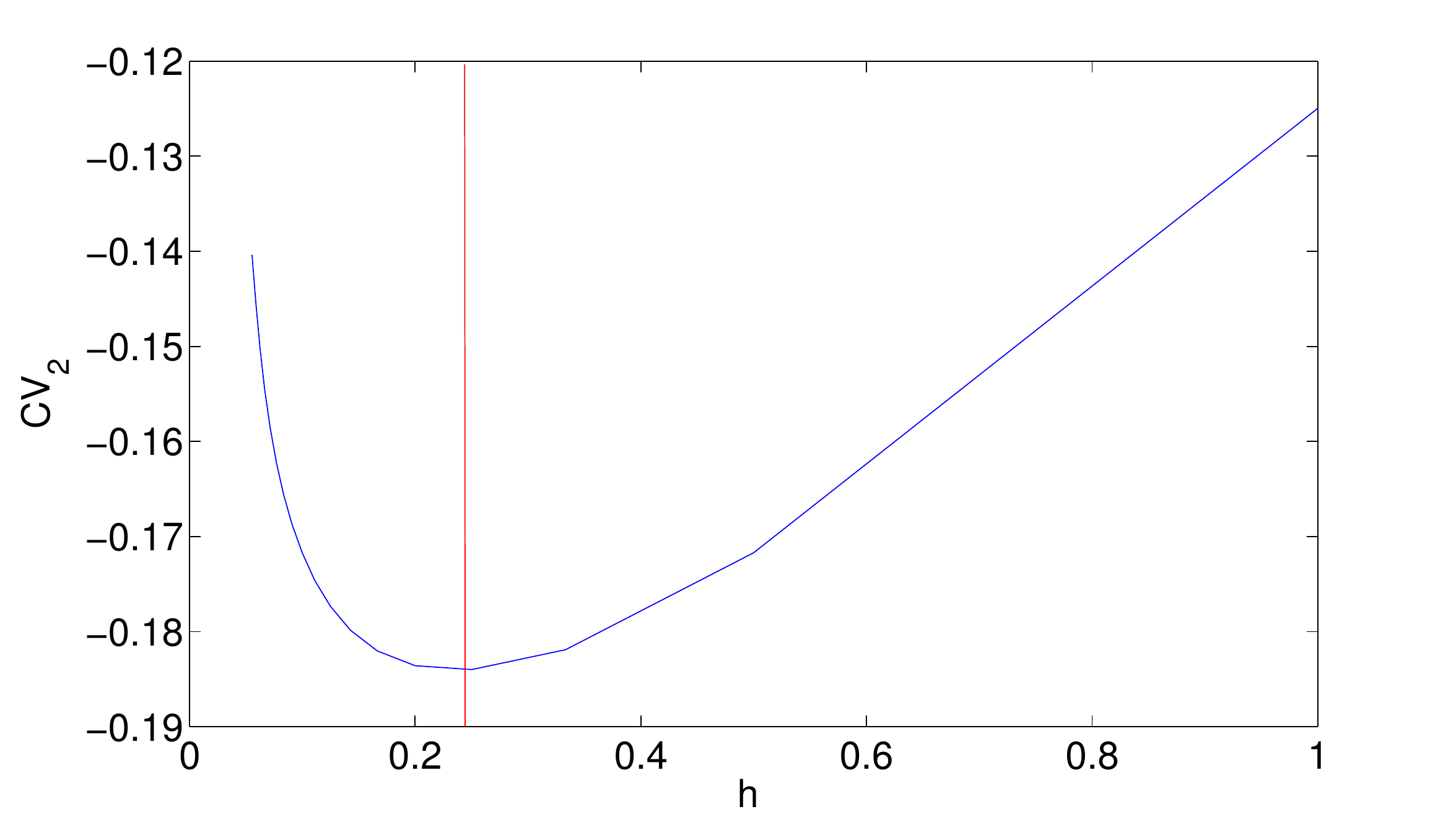}\\
a & b & c
\end{tabular}
\caption{Risks curves in function of $h$ for $f_{1,vM}$, $n=500$: a/ $R_{oracle}$  b/ $ R_{SPCO}$  c/ $CV_2$} Vertical red lines represent the bandwidth value $h$ minimizing each curve. 
 \label{courbes}
\end{center}\end{figure}

\begin{figure}[!h]\begin{center}
\begin{tabular}{ccc}
\includegraphics[scale=0.28]{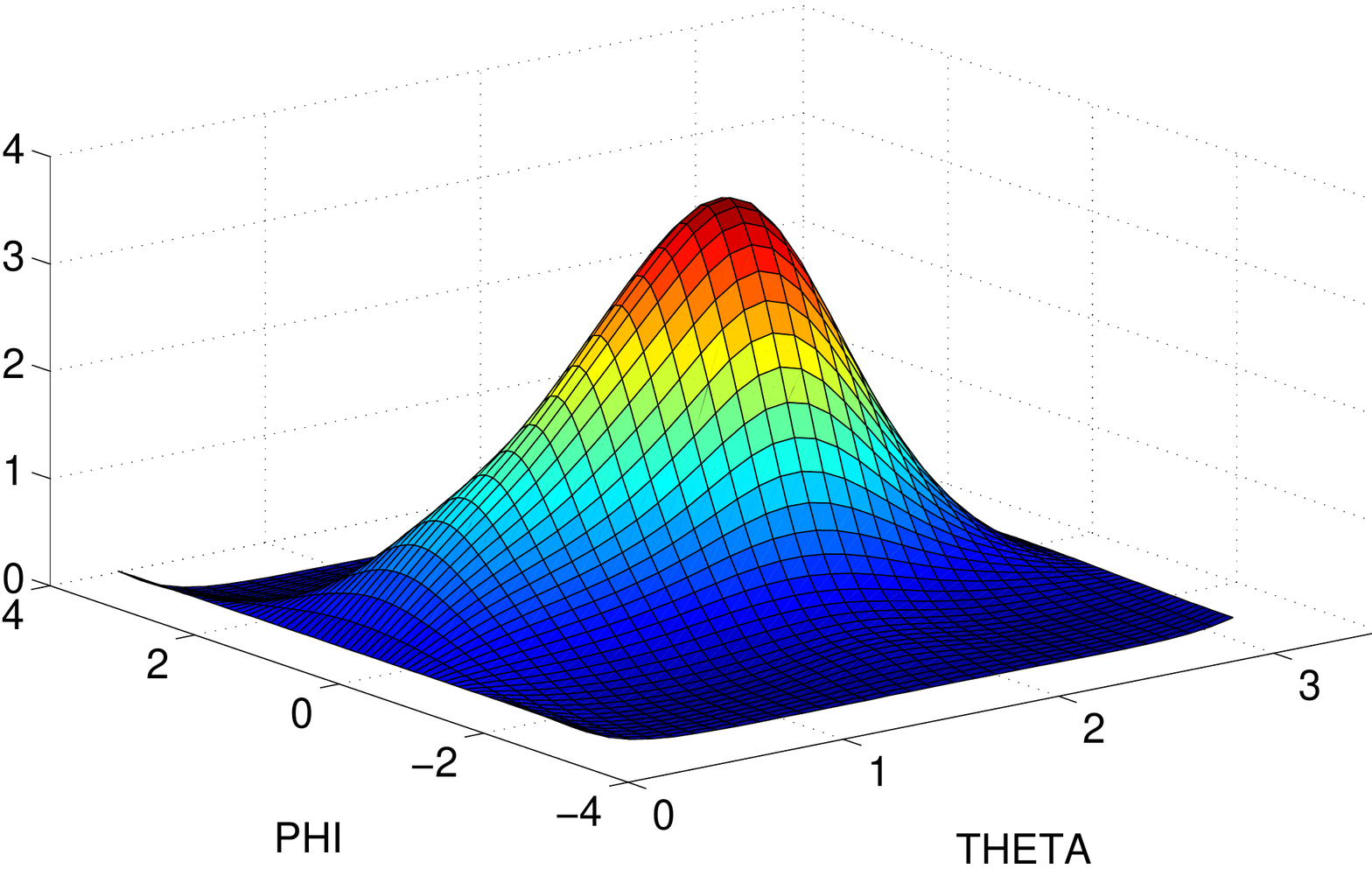} &
\includegraphics[scale=0.28]{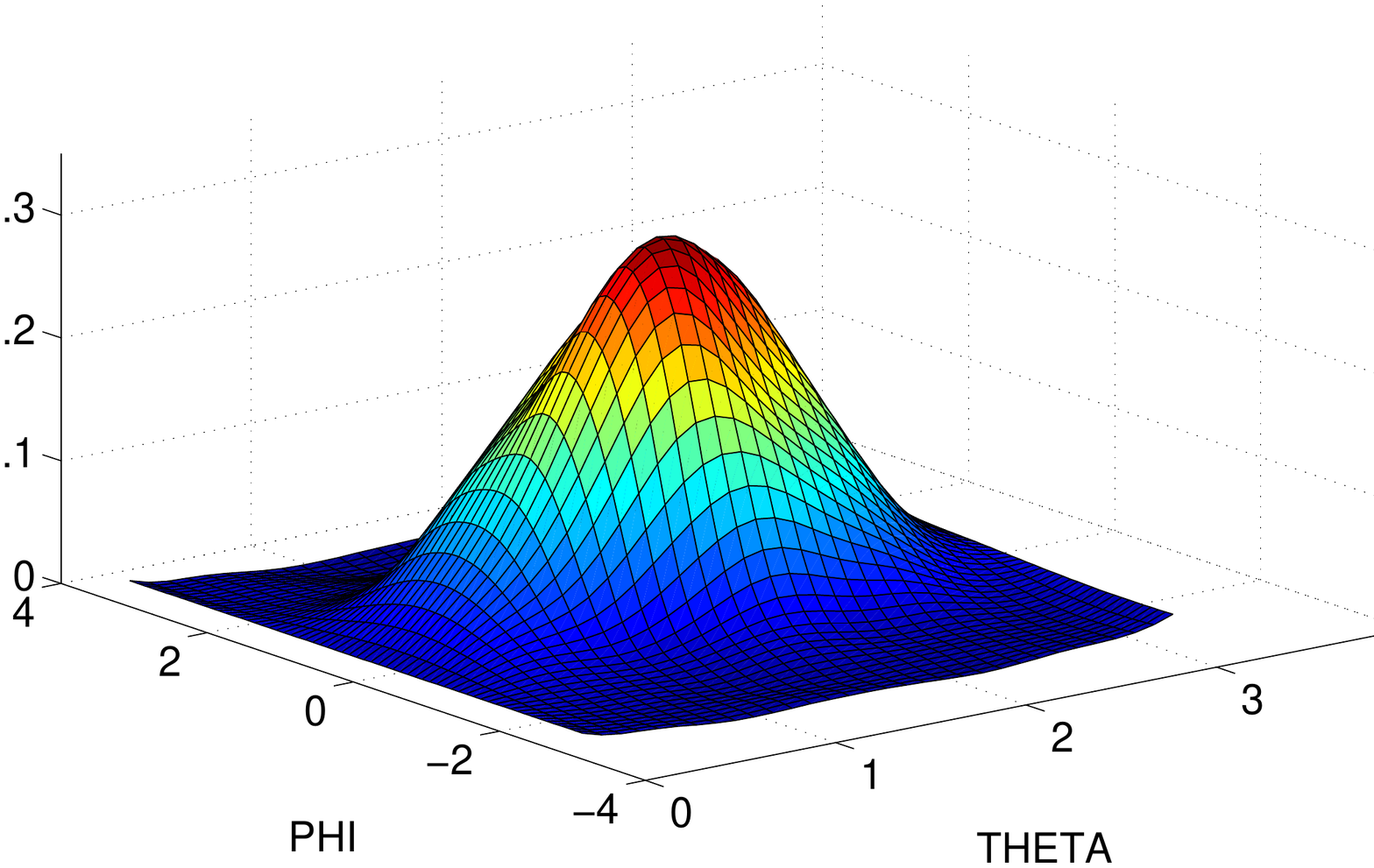} & 
\includegraphics[scale=0.28]{graphique-PCO.pdf} \\
a & b & c
\end{tabular}
\caption{  Reconstruction of $f_{1,vM}$, $n=500$:  a/  $\hat f_{h_{oracle}}$, $ h_{oracle}=0.33$  b/ SPCO  $\hat f_{\hat h}$, $\hat h=0.25$} c/ cross-validation $\hat f_{h_{CV2}}$,  $ h_{CV2}=0.25$
 \label{graphiques}
\end{center}\end{figure}

In the rest of this section, $SPCO$ will denote the estimation procedure related to $\hat f_{\check h}$. In Figure \ref{courbes}, for $n=500$ we plot as a function of $h$:   $ R_{oracle}:=\| \hat f_h -f_{1,vM}\|^2 -\| f_{1,vM}\|^2$ for the oracle,  $R_{SPCO}: =\| \hat f_h -\hat f_{\hmin}\|^2 + \pen(h) $ for SPCO and $CV_2(h)$ for cross-validation. We point out on each graphic the value of $h$ that minimizes each quantity. 
In Figure \ref{graphiques}, we plot in spherical coordinates, for $n=500$, the density $f_{1,vM}$ and density reconstructions for the oracle, SPCO and cross-validation. 
Eventually, in Tables \ref{MISE1} and \ref{MISE2}, we compute MISE to estimate $f_{1,vM}$ and $f_{2,vM}$ for  the oracle, SPCO and cross-validation  for $n=100$ and $n=500$, over $100$ Monte-Carlo runs. 
\bigskip

\begin{table}[!h]
\begin{center}
\begin{tabular}{| c | c | c | }
\hline
  & $n=100$ & $n=500$ \\
 \hline 
 Oracle & $0.0064$ & $0.0027$ \\
 \hline
 SPCO & 0.0091 &  0.0048 \\
 \hline 
 Cross -Validation & 0.0099 & 0.0053 \\
 \hline
\end{tabular}
\caption{ MISE  over 100 Monte-Carlo repetitions to estimate $f_{1,vM}$}
\label{MISE1}
\end{center}
\end{table}

\begin{table}[!h]
\begin{center}
\begin{tabular}{| c | c | c | }
\hline
  & $n=100$ & $n=500$ \\
 \hline 
 Oracle & $0.0051$ & $0.0027$ \\
 \hline
 SPCO & 0.0083 &  0.0043 \\
 \hline 
 Cross -Validation & 0.0096 & 0.0047 \\
 \hline
\end{tabular}
\caption{ MISE over 100 Monte-Carlo repetitions to estimate $f_{2,vM}$}
\label{MISE2}
\end{center}
\end{table}

When analyzing the results, SPCO shows quite satisfying performances. Indeed, SPCO is close to the oracle and is slightly better than cross-validation when looking at the MISE computations for both densities. 

\section{Proofs}\label{preuves}

Before proving Theorem \ref{theorem 2}, we need several intermediate results. The first one is the following proposition which is the counterpart of Proposition 4.1 of \cite{Lerasle} for $\S$.  

Let  
\begin{equation}\label{R1}
R_1(K):=2^{(d-3)/2} \sigma_{d-2}\int_{0}^{\infty} x^{(d-3)/2}K^2(x) dx.
\end{equation}

\begin{proposition}\label{lerasle2}
Assume that kernel $K$ satisfies Assumption 1. 
 Let $ \Upsilon \geq (1+2 \| f\|_{\infty}) \vee \frac{8 \pi \| K\|_{\infty}R_0(K)}{R_1(K)}$. There exists $n_0$, such that for $n \geq  n_0$ ($n_0$ not depending on $f$),  all $x \geq 1$ and for all $\eta \in (0,1)$ with probability larger than $1- \square |\mathcal{H}|  e^{-x}$, for all $h \in \mathcal{H}$ each of the following inequalities holds
\begin{eqnarray}
\| f -\hat f_h \|  &\leq& (1+\eta)  \left ( \|f-f_h \|^2  +\frac{c_0^2(h)c_2(h) }{n}\right) + \square \frac{\Upsilon x^2}{\eta^3 n} \label{lerasle21} \\
\| f-f_h\|^2 + \frac{ c_0^2(h)c_2(h)}{n} & \leq & (1+\eta) \| f -\hat f_h\|^2 + \square \frac{\Upsilon x^2}{\eta^3 n }. \label{lerasle22}
\end{eqnarray}
\end{proposition}

\textbf{Proof of Proposition \ref{lerasle2}}. 


To prove Proposition  \ref{lerasle2}, we need to verify  Assumptions (11) - (16) of \cite{Lerasle}. We remind that 
$$f_h= \E(\hat f_h) = c_0(h) \int_{y \in \S} f(y) K_{h^2}(x,y)   \omega_d(dy ).$$

Let us check Assumption (11) of \cite{Lerasle}. This one amounts to prove that for some $\Gamma$ and $\Upsilon$
$$
\Gamma (1+ \| f\|_{\infty})\vee \sup_{h \in \mathcal{H}} \| f_h \|^2 \leq \Upsilon.
$$

We have
\begin{eqnarray*}
\| f_h \|^2  & \leq  & \| f\|_{\infty} \int_x  \underbrace{   \left( \int_y c_0(h) K_{h^2}(x,y) \omega_d(dy) \right )    }_{=1}\left(  \int_y f(y) c_0(h) K_{h^2}(x,y) \omega_d(dy)\right)  \omega_d(dx)   \\
&\leq & \| f\|_{\infty} \int_y f(y)  \int_x c_0(h)  K_{h^2}(x,y)   \omega_d(dx)   \omega_d(dy)  \\
& \leq & \| f\|_{\infty},
\end{eqnarray*}
hence Assumption (11) in \cite{Lerasle} holds with $\Gamma=1$ and $\Upsilon \geq 1+ \| f\|_{\infty}$. 

Let us check Assumption (12) of \cite{Lerasle}. We have to prove that 
\[
 \int c_0^2(h)K^2_{h^2}(x,x) \omega_d(dx)  \leq \Upsilon n   \int \int c_0^2(h)K^2_{h^2}(x,y) \omega_d(dx) f(y) \omega_d(dy) .
\] 

 But since
$$  \int \int c_0^2(h)K^2_{h^2}(x,y) \omega_d(dx) f(y) \omega_d(dy) =  c^2_0(h)c_2(h) \int f(y) \omega_d(dy)= c^2_0(h)c_2(h),$$
and
$$
 \int c_0^2(h)K^2_{h^2}(x,x) \omega_d(dx) = 4\pi c_0^2(h) K^2(0),
$$
Assumption (12) amounts to check that 
$$
 4 \pi c_0^2(h) K^2(0) \leq \Upsilon n c^2_0(h)c_2(h) \Longleftrightarrow \Upsilon  \geq   \frac{4\pi K^2(0)}{n c_2(h)}.
$$
But using (\ref{c21}),
we have 
$$c_2(h)= R_1(K)h^{d-1}+o(1),$$
when $h$ tends to $0$ uniformly in $h$. Thus there exists $n_1$,
$ n_1$ independent of $f$, such that for $n\geq n_1$, $c_2(h) \geq \frac 12 R_1(K) h^{d-1}$.
Now using that $ h^{d-1} \geq \frac{\| K\|_{\infty}}{{R_0(K)n}}$ and $K(0) \leq \| K\|_{\infty}$, it is sufficient to have $\Upsilon \geq  \frac{8\pi\| K\|_{\infty}R_0(K)}{R_1(K)}$ to ensure Assumption (12) in \cite{Lerasle}.

Assumption (13) in \cite{Lerasle} consists to prove   that

$$ \|  f_h-f_{h'} \|_{\infty} \leq \Upsilon \vee \sqrt{\Upsilon n } \| f_h-f_{h'}\|_2.$$
For any $h\in \mathcal{H}$ and any $x \in \S$, we have
\[
\|f_h \|_{\infty}\leq  \| f \|_{\infty},
\]
therefore Assumption (13) in \cite{Lerasle} holds for 
$\Upsilon \geq 2 \| f \|_{\infty}$.

Assumptions (14) and (15) of \cite{Lerasle} consist in proving respectively that 
$$ \E \left [ c^2_0(h) \int  K_{h^2}(X,z) K_{h^2}(z,Y )\omega_d(dz)  \right ]^2 \leq \Upsilon c^2_0(h)c_2(h)  $$
and
$$
\sup_{x \in \S} \E \left [ c^2_0(h) \int  K_{h^2}(X,z) K_{h^2}(z,x) \omega_p(dz)  \right ]^2 \leq \Upsilon n. 
$$
We have 
\begin{eqnarray*}
 c_0^2(h) \int  K_{h^2}(x,z) K_{h^2}(z,y)  \omega_d(dz)   \leq  c_0^2(h)c_2(h) \wedge c_0(h) \| K\|_{\infty}.
\end{eqnarray*}
Indeed if $y=z$, then  $ \int  K_{h^2}(x,z) K_{h^2}(z,y)  \omega_d(dz) = \int  K^2_{h^2}(x,z)  \omega_d(dz) = c_2(h) $, otherwise
\[
 c_0^2(h) \int  K_{h^2}(x,z) K_{h^2}(z,y)  \omega_d(dz) \leq c_0(h) \| K\|_{\infty}\underbrace{ c_0(h) \int  K_{h^2}(x,z)   \omega_d(dz)}_{=1}= c_0(h) \| K\|_{\infty}.
\]
Furthermore (\ref{c01}) entails that there exists $n_2$ independent of $f$, such that for $n \geq n_2$, $c^{-1}_0(h) \geq \frac 1 2 R_0(K) h^{d-1} $ and consequently $c_0(h) \leq \frac{2n}{\| K\|_{\infty}}$, using (\ref{H}). Thus for  $n \geq n_2$
\begin{equation}\label{eq3}
 c_0^2(h) \int  K_{h^2}(x,z) K_{h^2}(z,y)  \omega_d(dz)   \leq c_0^2(h)c_2(h) \wedge 2n. 
\end{equation}

We have
\begin{eqnarray}
\E  \left [  c^2_0(h) \int_z  K_{h^2}(X,z) K_{h^2}(z,x) \omega_d(dz) \right ] &=&  c_0^2(h) \int_z \left( \int_{y}   K_{h^2}(y,z)f(y)  \omega_d(dy) \right) K_{h^2}(z,x)  \omega_d(dz) \nonumber \\
&\leq & \| f\|_{\infty }  c_0(h)  \int_z    \underbrace{c_0(h) \int_{y}   K_{h^2}(y,z)  \omega_d(dy) }_{=1}K_{h^2}(z,x)  \omega_d(dz) \nonumber  \\
&\leq &  \| f\|_{\infty} \label{eq4}.
\end{eqnarray}
 
Therefore for $n \geq n_2,$
 
 \begin{eqnarray*}
\sup_{x \in \S} \E \left [ c^2_0(h) \int  K_{h^2}(X,z) K_{h^2}(z,x) \omega_d(dz)  \right ]^2 &\leq&  \sup_{(x,y)}    \left ( c^2_0(h) \int  K_{h^2}(x,z) K_{h^2}(z,y) \omega_d(dz)   \right ) \\ 
& \quad \times & \sup_{x} \E \left  [ c^2_0(h) \int  K_{h^2}(X,z) K_{h^2}(z,x) \omega_d(dz) \right  ]\\
&\leq & (c_0^2(h)c_2(h) \wedge 2n)  \| f\|_{\infty},
\end{eqnarray*} 
using (\ref{eq3}) and (\ref{eq4}).

And we have
 \begin{eqnarray*}
\E \left [ c^2_0(h) \int  K_{h^2}(X,z) K_{h^2}(z,Y )\omega_d(dz)  \right ]^2 & \leq &\sup_{x} \E \left  [  c^2_0(h) \int  K_{h^2}(X,z) K_{h^2}(z,x) \omega_p(dz)  \right ]^2 \\
& \leq & (c_0^2(h)c_2(h) \wedge 2n)  \| f\|_{\infty}
\end{eqnarray*} 
using (\ref{eq3}) and (\ref{eq4}).

Hence Assumption  (14) and (15) in \cite{Lerasle}  hold for $\Upsilon \geq 2\| f\|_{\infty}$. 

Now  let 
$t \in \mathbb{B}_{{c_0(h)K_{h^2}}}$ is the set of functions $t$  which can be written $t(x)= \int a(z) c_0(h) K_{h^2}(z,x) \omega_d(dz)$  for some $a \in \mathbb{L}^2(\S)$ with $\| a\| \leq 1 $. Now let $a  \in L^2(\S)$ be such that $\| a\|=1$ and $t(y)= \int a(x) c_0(h) K_{h^2}(x,y) \omega_d(dx)$ for all $y \in \S$. 
To verify Assumption (16) in \cite{Lerasle} we have to prove that 
$$ \sup_{t \in  \mathbb{B}_{{c_0(h)K_{h^2}}} } \int t(x) f(x) \omega_d(dx) \leq \Upsilon \vee \sqrt{\Upsilon c_0^2(h) c_2(h)}.$$
Using Cauchy Scharwz inequality one gets 
\[
t(y) \leq  \sqrt{ \int_{\S} a^2(x) \omega_d(dx)} \sqrt{ c_0^2(h) \int_{\S} K^2_{h^2}(x,y) \omega_d(dx) }  \leq  \sqrt{c_0^2(h) c_2(h)}.
\]
Thus for any $t \in  \mathbb{B}_{{c_0(h)K_{h^2}}} $
\begin{eqnarray}
 \int t^2(x) f(x) \omega_d(dx)  \leq \| t\|_{\infty} \langle |t|,f \rangle \leq  \sqrt{c_0^2(h) c_2(h)}  \| f\| \| t\|,
 \end{eqnarray}
but using Cauchy Schwarz inequality and Fubini, one gets 

\begin{eqnarray*}
\| t\| &= &\int_x \left (\int_y a(y)c_0(h) K_{h^2}(x,y) \omega_d(dy) \right )^2   \omega_d(dx) \\
& \leq & \int_x \left(  \int_y a^2(y) c_0(h) K_{h^2}(x,y) \omega_d(dy) \right )\left( \int_y  c_0(h) K_{h^2}(x,y) \omega_d(dy) \right)  \omega_d(dx) \\
&= & \int_x \int_y a^2(y) c_0(h) K_{h^2}(x,y) \omega_d(dy)  \omega_d(dx) = \int_y a^2(y) \left( \int_x  c_0(h) K_{h^2}(x,y) \omega_d(dx)  \right) \omega_d(dy)   =1.
\end{eqnarray*}

And 

\begin{eqnarray*}
 \int t^2(x) f(x) \omega_p(dx)  &\leq & \| f\| \sqrt{c_0^2(h) c_2(h)}   \\
&\leq &  \sqrt{\| f\|_{\infty}} \sqrt{c_0^2(h) c_2(h)} \\
& \leq & \sqrt{\Upsilon c_0^2(h) c_2(h)},
 \end{eqnarray*}
 hence Assumption (16) in \cite{Lerasle} is verified. 

Finally, assumptions (11)-(16) from \cite{Lerasle} hold in the spherical setting, for $n\geq n_0= \max(n_1, n_2)$ and  if  $\Gamma=1$ and
\[
\Upsilon \geq (1+2 \| f\|_{\infty}) \vee \frac{8 \pi \| K\|_{\infty}R_0(K)}{R_1(K)}.
\]
This enables us to use Proposition 4.1 of \cite{Lerasle} which gives Proposition \ref{lerasle2}. This ends the proof of Proposition \ref{lerasle2}. 

\vspace{5mm}

%


The next proposition is the second step to prove Theorem \ref{theorem 2}.
\begin{proposition}\label{theorem9}
Assume that the kernel $K$ satisfies Assumption 1 and $\| f\|_{\infty} < +\infty$. Let $x \geq 1$ and $\theta \in (0,1)$. With probability larger than $1-C_1 |\mathcal{H}| e^{-x}$ , for any $h \in \mathcal{H}$,

\begin{equation}
\begin{split}
(1-\theta) \| \hatfhat- f \|^2  &\leq (1+\theta) \| \hat f_h -f\|^2 + \left( pen_{\lambda}(h) -2 \frac{\langle c_0(h) K_{h^2},c_0(\hmin) K_{\hmin^2} \rangle}{n} \right )  \\
& \quad \quad  -\left( pen_{\lambda}(\hat h) - 2 \frac{\langle c_0(\hat h) K_{\hat h^2},c_0(\hmin) K_{\hmin^2} \rangle}{n} \right ) \\
& \quad \quad + \frac{C_2}{\theta}\| f_{\hmin} -f\|^2 + \frac{C(K)}{\theta} \left(  \frac{\| f\|_{\infty} x^2}{n} + \frac{x^3c_0(\hmin)}{n^2}\right),
\end{split}
\end{equation}
where $C_1$ and $C_2$ are absolute constants and $C(K)$ only depends on $K$. 
\end{proposition}

In order to avoid any confusion, we recall that $K_{h^2}=K_{h^2}(\cdot , \cdot )$ and 
$$
\langle c_0(h) K_{h^2},c_0(\hmin) K_{\hmin^2} \rangle = \int_{\S} c_0(h) K_{h^2}(x,y)c_0(\hmin) K_{\hmin^2}(x,y) \omega_d(dy). 
$$
Once again, we would like to draw the attention that the quantity $\int_{\S} c_0(h) K_{h^2}(x,y)c_0(\hmin) K_{\hmin^2}(x,y) \omega_d(dy)$ does not depend on $x$. Indeed, we have, using Remark \ref{remarque1}
\begin{eqnarray*}
\int_{\S} K_{h^2}(x,y) K_{\hmin^2}(x,y) \omega_d(dy) &= &   \int_{\S}  K \left (\frac{1-x^Ty}{h^2} \right ) K \left (\frac{1-x^Ty}{\hmin^2} \right ) \omega_d(dy) \\
&=&  \sigma_{d-2} \int_{-1}^{1} K \left (\frac{1-t}{h^2} \right) K \left (\frac{1-t}{\hmin^2} \right ) (1-t^2)^{(d-3)/2}dt . 
\end{eqnarray*}

\textbf{Proof of Proposition \ref{theorem9}.}  The proof  follows the proof of Theorem 9 in \cite{Lacour} adapted to $\S$. 
Let $\theta' \in(0,1)$ be fixed and chosen later. Using the definition of $\hat h$, we can write, for any $h \in \mathcal{H}$
\begin{eqnarray*} 
\| \hatfhat -f\|^2 +pen_\lambda(\hat h) 
&= &  \| \hatfhat -\hat f_{\hmin}\|^2 + \mbox{pen}_\lambda(\hat h) + \| \hat f_{\hmin} - f\|^2  +2 \langle \hatfhat -\hat f_{\hmin}, \hat f_{\hmin} - f\rangle  \\
& \leq &  \| \hat f_h -\hat f_{\hmin}\|^2 + pen_\lambda(h) + \| \hat f_{\hmin} - f\|^2  +2 \langle \hatfhat -\hat f_{\hmin}, \hat f_{\hmin} - f\rangle \\
& \leq & \| \hat f_h -f\|^2 +2 \|  f- \hat f_{\hmin} \|^2 +2 \langle \hat f_h -f, f-\hat f_{\hmin} \rangle + pen_\lambda(h) + 2
\langle \hatfhat -\hat f_{\hmin}, \hat f_{\hmin}- f \rangle
\end{eqnarray*}

Consequently 
\begin{equation}\label{16}
\begin{split}
\| \hatfhat -f\|^2 & \leq \| \hat f_h -f\|^2 +\left(   pen_\lambda(h)  -2 \langle \hat f_h -f,  \hat f_{\hmin}-f \rangle  \right)  \\
&\quad \quad -  \left(  pen_\lambda(\hat h) - 2  \langle \hatfhat -f, \hat f_{\hmin}- f \rangle \right).
\end{split}
\end{equation}

Then for a given $h$, we study the term 
$$  2 \langle \hat f_h -f,  \hat f_{\hmin}-f \rangle.  $$

Let us introduce the degenerate $U$-statistic
$$U(h,\hmin) := \sum_{i \neq j } \left \langle c_0(h)K_{h^2}(. \; , X_i) -f_h, c_0(\hmin) K_{\hmin}(. \;, X_j) -f_{\hmin} \right \rangle $$

and the following centered variable
$$ V(h,h'):= \langle \hat f_h -f_h, f_{h'}-f \rangle.  $$

We first center the terms  
\begin{eqnarray*}
\langle \hat f_h -f,  \hat f_{\hmin}-f \rangle &=& \langle \hat f_h -f_h, \hat f_{\hmin} -f_{\hmin} \rangle + V(h, \hmin) + V(\hmin,h) + \langle f_h -f, f_{\hmin} -f\rangle.
\end{eqnarray*}

Now 
\begin{eqnarray*}
\langle \hat f_h -f_h, \hat f_{\hmin}- f_{\hmin} \rangle 
&=&   \frac{1}{n^2} \sum_{i,j=1}^n \langle c_0(h) K_{h^2}(. \;, X_i) -f_h, c_0(\hmin)K_{\hmin}(., X_j) -f_{\hmin} \rangle   \\
&=&  \frac{1}{n^2} \sum_{i=1}^n \langle c_0(h) K_{h^2}(. \;, X_i) -f_h, c_0(\hmin) K_{\hmin}(., X_i) -f_{\hmin} \rangle+\frac{U(h,\hmin)}{n^2}.
\end{eqnarray*}

Then 

\begin{eqnarray*}
\langle \hat f_h -f_h, \hat f_{\hmin}-f_{\hmin} \rangle =  \frac{  \langle c_0(h)K_h^2, c_0(\hmin) K_{\hmin^2}\rangle}{n}   -\frac 1 n \langle  \hat f_h, f_{\hmin}\rangle -\frac 1 n \langle f_h, \hat f_{\hmin} \rangle +\frac 1 n \langle f_h, f_{\hmin} \rangle +\frac{U(h,\hmin)}{n^2}.
\end{eqnarray*}

Finally we obtain 
\begin{eqnarray}
\langle \hat f_h -f, \hat f_{\hmin}-f \rangle &=&  \frac{  \langle c_0(h)K_h^2, c_0(\hmin) K_{\hmin^2}\rangle}{n} +\frac{U(h,\hmin)}{n^2} \label{eq0} \\
 &-& \frac 1 n \langle  \hat f_h, f_{\hmin}\rangle -  \frac 1 n \langle   f_h, \hat f_{\hmin}\rangle  + \frac{1}{n} \langle f_h, f_{\hmin}\rangle \label{eq1} \\ 
 &+& V(h,\hmin) +V(\hmin,h) + \langle f_h -f, f_{\hmin}-f \rangle. \label{eq2}
\end{eqnarray}

We first control the last term of (\ref{eq0}) involving a U-statistics. This is done in the next lemma. 
\begin{lemma} \label{Ustat}
 With probability greater than $1-5.54 |\mathcal{H}| e^{-x} $, for any $h$ in $ \mathcal{H}$,
 \[
\left | \frac{U(h,\hmin)}{n^2}\right  |\leq   \frac{ \theta' c_0^2(h) c_2(h) }{n} + \frac{\square \| f\|_{\infty}x^2}{\theta' n} + \frac{\square c_0(\hmin)  \| K \|_{\infty}x^3} {\theta' n^2}.
\]

\end{lemma}

\textbf{Proof of Lemma \ref{Ustat}.}

We have 

\begin{eqnarray*}
 U(h,\hmin) & = &\sum_{i \neq j} \langle c_0(h) K_h^2(., X_i) -f_h, c_0(\hmin) K_{\hmin^2}(., X_j) - f_{\hmin}  \rangle  \\
 &=& \sum_{i=2}^n \sum_{j<i} G_{h,\hmin} (X_i,X_j) +G_{\hmin, h} (X_i,X_j)  
 \end{eqnarray*}
where
$$ G_{h,h'}(s,t) =\langle c_0(h) K_{h^2}(.,s)-f_h, c_0(h') K_{h'^2}(.,t) -f_{h'}\rangle . $$
We apply Theorem 3.4 of \cite{Houdre}
$$ \P \left  ( | U(h,\hmin) | \geq \square \left ( C\sqrt{x} +Dx +B x^{3/2} +Ax^2\right) \right) \leq 5.54 e^{-x},$$
with $A,B,C$ and D defined subsequently. First, we have
\begin{eqnarray*}
 \| f_{\hmin}  \|_{\infty}&=& \| \E(\hat f_{\hmin})  \|_{\infty }  = \| c_0(\hmin) \int_{\S} K_{\hmin^2}(x,y)f(y) \omega_d(dy) \|_{\infty} \\
 &\leq &  c_0(\hmin) \| K\|_{\infty} \int_{\S} f(y) \omega_d(dy) \\
 &\leq &   c_0(\hmin)   \| K\|_{\infty}. \\
 \end{eqnarray*}
 We have
 \begin{eqnarray*}
 A&:=& \|  G_{h,\hmin} +G_{\hmin,h}\|_{\infty} \\
 &\leq & \| G_{h,\hmin}\|_{\infty} + \| G_{\hmin,h}\|_{\infty} = 2 \| G_{h,\hmin}\|_{\infty},
 \end{eqnarray*}
because $  G_{\hmin,h}=  G_{h,\hmin} $.
 We have
 \begin{equation}
 \begin{split}
  \| G_{h,\hmin} \|_{\infty} & =  \sup_{s,t} \left | \int_{\S} \left  ( c_0(h)K_{h^2}(u,s)-f_h(u) \right  ) \left (c_0(\hmin)K_{\hmin^2}(u,t) -f_{\hmin}(u) \right) \omega_d(du) \right | \\
  &\leq  \sup_{u,t} | c_0(\hmin) K_{\hmin^2}(u,t)-f_{\hmin}(u) | \sup_{s} \int |c_0(h)K_{h^2}(u,s) -f_h(u)| \omega_d(du)  \\
  & \leq  \left ( c_0(\hmin) \|  K \|_{\infty} + \|  f_{\hmin} \|_{\infty} \right )  \Big ( \sup_{s} c_0(h)  \int K_{h^2}(u,s) \omega_d(du) \\ 
  & \quad \quad \quad + c_0(h)  \int \int K_{h^2}(u,y)f(y) \omega_d(dy) \omega_d(du) \Big )  \\
  &\leq   2 c_0(\hmin)  \| K \| _{\infty}  \left (1 +  \int f(y) \underbrace{c_0(h) \int K_{h^2}(u,y)  \omega_d(du)}_{=1} \omega_d(dy) \right ) \\
  & \leq  4 c_0(\hmin)  \| K \| _{\infty} . \\
   \end{split}
 \end{equation}
Consequently we have  that 
 $$ A \leq 8 c_0(\hmin)  \| K\|_{\infty }$$
and
$$ \frac{Ax^2}{n^2} \leq \frac{8x^2 c_0(\hmin) \| K \|_{\infty} }{n^2}.$$
 We define 
 $$B^2 := (n-1) \sup_t \E [ (G_{h,\hmin}(t,X_2) +G_{\hmin,h}(t,X_2))^2].$$
For any $t$ we have
\begin{equation}
\begin{split}
\E [G^2_{h,\hmin}(t,X_2)] &= \E \left [ \left (\int (c_0(h) K_{h^2}(u,t) -f_h(u))(c_0(\hmin) K_{\hmin^2} (u,X_2)  -\E [c_0(\hmin) K_{\hmin^2}(u,X_2)] )  \omega_d(du)\right )^2 \right ] \\
& \leq  \E  \bigg [ \int (c_0(h)K_{h^2}(u,t) -f_h(u))^2 \omega_d(du) \times  \\
& \quad \int  (c_0(\hmin)K_{\hmin^2}(u,X_2) -\E[c_0(\hmin) K_{\hmin^2}(u,X_2)])^2  \omega_d(du) \bigg  ]  \\
& \leq  2  \bigg (\int c_0^2(h)K_{h^2}^2(u,t) \omega_d(du)  +  \int f_h^2(u) \omega_d(du) \bigg) \times  \\
& \quad \int  \E  \bigg [c_0(\hmin)K_{\hmin^2}(u,X_2) -\E[c_0(\hmin) K_{\hmin^2}(u,X_2)] \bigg ]^2  \omega_d(du) \\
&\leq  2 \times  \bigg( \int c_0^2(h) K^2_{h^2}(u,t)  \omega_d(du) + \int_{u} \big (c_0(h)\int_y  K_{h^2}(u,y)f(y) \omega_d(dy )\big  )^2 \; \;  \omega_d(du) \bigg) \\
& \quad \times \int \E [ c_0^2(\hmin)K^2_{\hmin}(u,X_2)] \omega_d(du)   \\
& \leq 2 \left  ( c_0^2(h) c_2(h) +  \int \int c_0^2(h) K^2_{h^2}(u,y) f(y)\omega_d(dy) \int f(y) \omega_d(dy) \omega_d(du) \right ) \\
& \quad \times \int \E [ c_0^2(\hmin)K^2_{\hmin^2}(u,X_2)] \omega_d(du) \\
& \leq  4 c_0^2(h) c_2(h)c^2_0(\hmin) c_2(\hmin).
\end{split}
\end{equation}

Therefore 

$$ B^2 \leq 8 (n-1) c_0^2(h) c_2(h)c^2_0(\hmin) c_2(\hmin)$$

and 

$$\frac{B^2x^3}{n^4} \leq \frac{8}{n} c_0^2(h) c_2(h)  \frac{1}{n^2} c^2_0(\hmin) c_2(\hmin) x^3 .$$

Now using 
$$\sqrt{ab} \leq \theta \frac a 2  + \theta^{-1} \frac b 2,$$

we obtain 

$$\frac{Bx^{3/2}}{n^2} \leq  \frac{\theta'}{3n} c_0^2(h) c_2(h) +  \frac{6}{\theta'} c^2_0(\hmin) c_2(\hmin) \frac{x^3}{n^2}.$$

Now we have

\begin{equation*}
\begin{split}
C^2 &:= \sum_{i=2}^n \sum_{j=1}^{i-1} \E [(G_{h,\hmin})(X_i,X_j) +G_{\hmin,h}(X_i,X_j))^2 ] \\
& \leq \square \times n^2 \E[G^2_{h,\hmin}(X_1,X_2)] \\
& = \square n^2 \E \left [ \left  (  \int (c_0(h)K_{h^2}(u,X_1)-f_h(u))(c_0(\hmin)K_{\hmin^2}(u,X_2)-f_{\hmin}(u)) \omega_d(du) \right )^2 \right] \\
&= \square \times n^2 \E  \Bigg [  \Bigg ( \int c_0(h) K_{h^2}(u,X_1) c_0(\hmin)K_{\hmin^2}(u,X_2) \omega_d(du)   \\
& \qquad \qquad \qquad  - \int c_0(h)K_{h^2}(u,X_1) \left( \int c_0(\hmin) K_{\hmin^2}(u,y) f(y) \omega_d(dy) \right)  \omega_d(du) \\
& \qquad \qquad \qquad   - \int c_0(\hmin) K_{\hmin^2}(u,X_2)   \left(\int_{} c_0(h) K_{h^2}(u,y)f(y) \omega_d(dy )  \right) \omega_d(du) \\
&  \qquad \qquad \qquad  + \int_u  \left( \int_{y} c_0(\hmin) K_{\hmin^2}(u,y)f(y) \omega_d(dy) \right)  \left (\int_{y} c_0(h) K_{h^2}(u,y)f(y) \omega_d(dy ) \right )\omega_d(du)  \Bigg )^2\Bigg ].\\
& \leq    \square n^2(A_1+A_2+A_3+A_4).
\end{split}
\end{equation*}

We have for $A_2$
\begin{equation*}
\begin{split}
 & \E  \left [ \int c_0(h)K_{h^2}(u,X_1) \left( \int c_0(\hmin) K_{\hmin^2}(u,y) f(y) \omega_d(dy) \right)   \omega_d(du) \right ]^2   \\
&  \leq \| f \|_{\infty}^2  \E \left ( \int c_0(h) K_{h^2}(u,X_1) \underbrace{ \left( \int c_0(\hmin) K_{\hmin^2}(u,y) \omega_d(dy) \right)}_{=1}  \omega_d(du) \right )^2 \\
& =  \| f \|_{\infty}^2 \E \left ( \int c_0(h) K_{h^2}(u,X_1) \omega_d(du) \right )^2 \\
& \leq  \| f \|_{\infty}^2 \int \left (  \int c_0(h) K_{h^2}(u,y) \omega_d(du) \right )^2 f(y)  \omega_d(dy) =  \| f \|_{\infty}^2.
\end{split}
\end{equation*}
With similar computations, we obtain the same bound for $A_3$. As for $A_4$, we get
\begin{equation*}
\begin{split}
 &  \E \left [   \int_u  \left(  \int c_0(\hmin) K_{\hmin^2}(u,y)f(y) \omega_d(dy ) \right) \left( \int c_0(h) K_{h^2}(u,y)f(y) \omega_d(dy )  \right)\omega_d(du)   \right ]^2 \\
& \leq  \E  \left [ \int \| f\|_{\infty}  \underbrace{ \left(\int c_0(\hmin) K_{\hmin^2}(u,y) \omega_d(dy) \right)}_{=1} \left( \int c_0(h) K_{h^2}(u,y)f(y) \omega_d(dy ) \right) \omega_d(du)    \right ]^2   \\
& \leq \| f\|_{\infty}^2  \E \left [ \int \int c_0(h)K_{h^2}(u,y)f(y) \omega_d(dy ) \omega_d(du )  \right ]^2 \\
& \leq   \| f\|_{\infty}^2 \left (  \int f(y) \int c_0(h)K_{h^2}(u,y) \omega_d(du )  \omega_d(dy ) \right )^2 =  \| f\|_{\infty}^2. 
\end{split}
\end{equation*}

Hence 

\begin{equation}
C^2 \leq \square \times n^2 \E \left [ \left (  \int c_0(h) K_{h^2}(u,X_1)c_0(\hmin)K_{\hmin}(u,X_2) \omega_d(du)\right)^2   \right ] + \square  \| f\|^2_{\infty} \times n^2.
\end{equation}

It remains to bound $A_1$. We have

\begin{equation*}
\begin{split}
\E  &\left [ \left (  \int c_0(h) K_{h^2}(u,X_1)c_0(\hmin)K_{\hmin}(u,X_2) \omega_d(du)\right)^2   \right ] \\
&= \int_y \int_x \left ( \int_u c_0(h)K_{h^2}(u,x)c_0(\hmin)K_{h^2}(u,y) \omega_d(du) \right  )^2 f(x)\omega_d(dx) f(y) \omega_d(dy)  \\
& \leq \| f\|_{\infty} \int_y \int_x \left (  \int_u c_0^2(h) K^2_{h^2}(u,x)c_0(\hmin) K_{\hmin^2}(u,y) \omega_d(du) \int_u c_0(\hmin)  K_{\hmin^2}(u,y) \omega_d(du)\right )\omega_d(dx) f(y) \omega_d(dy) \\
& \leq  \| f\|_{\infty} \int_{y} \int_x \int_u c^2_0(h) K^2_{h^2}(u,x)c_0(\hmin) K_{\hmin^2}(u,y) \omega_d(du) \omega_d(dx) f(y) \omega_d(dy) \\
& \leq  \| f\|_{\infty} c_0^2(h) \int_y \int_u \int_x K^2_{h^2}(u,x) \omega_d(dx) c_0(\hmin) K_{\hmin^2}(u,y) f(y) \omega_d(dy) \omega_d(du) \\
& \leq   \| f\|_{\infty} c_0^2(h) c_2(h) \int_y f(y)  \underbrace{\left( \int_u  c_0(\hmin) K_{\hmin^2}(u,y) \omega_d(du)  \right )}_{=1} \omega_d(dy)  \\
& \leq  \| f\|_{\infty} c_0^2(h) c_2(h). 
\end{split}
\end{equation*}

Finally,
$$
C \leq \square \times n \| f\|_{\infty}^{\frac 1 2} c_0(h) \sqrt{c_2(h)} +\square \| f\|_{\infty} \times n,
$$

hence, since $x \geq 1 $ we get 

\begin{eqnarray*}
\frac{C\sqrt{x}}{n^2} &\leq&  \frac{\square  \| f\|_{\infty}^{\frac 1 2} c_0(h) \sqrt{c_2(h)} \sqrt{x}}{n} +\frac{ \square \| f\|_{\infty} \sqrt{x}}{n} \\
&\leq & \frac{\theta'   c_0^2(h) {c_2(h)}}{3n} +\frac{ \square  \| f\|_{\infty} x}{\theta' n} + \frac{ \square \| f\|_{\infty} \sqrt{x}}{ n} \\
& \leq & \frac{\theta'  c_0^2(h) {c_2(h)}}{3n} +\frac{ \square  \| f\|_{\infty} x}{\theta' n}.
\end{eqnarray*}



Now let us consider 

\[ 
\mathcal{S} := \left \{ a=(a_i)_{2\leq i \leq n}, b=(b_i)_{1\leq i \leq n-1} : \sum_{i=2}^n\ \E[a_i^2(X_i)] \leq 1 , \sum_{i=1}^{n-1} \E [b_i^2(X_i)] \leq 1  \right \}.
\]

We have
\[
D:= \sup_{(a,b) \in \mathcal{S}} \left \{ \sum_{i=2}^n \sum_{j=1}^{i-1} \E [ ( G_{h,\hmin}(X_i, X_j) + G_{\hmin, h}(X_i,X_j)) a_i(X_i) b_j(X_j) ] \right \}.
\]

We have for $(a,b) \in \mathcal{S}$
\begin{equation}
\begin{split}
\sum_{i=2}^n  & \sum_{j=1}^{i-1}  \E[G_{h, \hmin} (X_i,X_j) a_i(X_i) b_j(X_j)] \\
 & \leq  \sum_{i=2}^n   \sum_{j=1}^{i-1}  \E \int | c_0(h) K_{h^2}(u,X_i) - f_h(u) | |a_i(X_i)|  \\
 & \quad \qquad \qquad   \times  |c_0(\hmin)K(u,X_j) - f_{\hmin}(u) | |b_{j}(X_j)| \omega_d(du) \\
 & \leq \sum_{i=2}^n \sum_{j=1}^{n-1} \int \E [ | c_0(h) K_{h^2}(u,X_i) - f_h(u) | |a_i(X_i)|  ] \\
 & \quad \qquad \qquad   \times   \E [ |c_0(\hmin)K_{\hmin^2}(u,X_j) - f_{\hmin}(u) | |b_{j}(X_j)| ]  \omega_d(du)  \
\end{split}
\end{equation}
and for any $u$, using  Cauchy Schwarz inequality

\begin{equation*}
\begin{split}
& \sum_{i=2}^n \E |c_0(h) K_{h^2}(u,X_i) - f_h(u)| |a_i(X_i)|  \\
& \leq \sqrt{n} \sqrt{\sum_{i=2}^n \E^2 [ c_0(h) K_{h^2}(u,X_i) -f_h(u)| |a_i(X_i)|] } \\
& \leq \sqrt{n} \sqrt{  \sum_{i=2}^n \E [ |c_0(h) K_{h^2}(u,X_i) - f_h(u) |^2] \E[a_i^2(X_i)]}  \\
& \leq \sqrt{n}  \sqrt{ \sum_{i=2}^n \E [c_0^2(h) K_{h^2}^2(u,X_i)] \E[a_i^2(X_i)]  } \\
& \leq  \sqrt{n} \sqrt{ \| f\|_{\infty} c_0^2(h)c_2(h) \sum_{i=2}^n  \E[a_i^2(X_i)] } \\
& \leq \sqrt{c_0^2(h)c_2(h) } \sqrt{n \| f\|_{\infty}}.
\end{split}
\end{equation*}
Now since
\begin{equation*}
\int f_{\hmin} (u) \omega_d(du)=1,
\end{equation*}
and 
\[
\int  \E  [ c_0(\hmin) K_{\hmin^2} (u,X_j)] \omega_d(du) =1,
\]
we have 
\begin{equation*}
\begin{split}
& \sum_{j=1}^{n-1} \int \E \left [ | c_0(\hmin) K_{\hmin^2} (u,X_j) -f_{\hmin} (u) |  \omega_d(du) |b_j(X_j)| \right ] \\
&\leq  2   \sum_{j=1}^{n-1} \E [|b_j(X_j)|] \\
& \leq  2   \sqrt{n} \sqrt{ \sum_{j=1}^{n-1} \E [|b^2_j(X_j)|]} \\
& \leq 2 \sqrt{n}.
\end{split}
\end{equation*}

Finally,
\[
\sum_{i=2}^n  \sum_{j=1}^{i-1}  \E[G_{h, \hmin} (X_i,X_j) a_i(X_i) b_j(X_j)]  \leq 2 n \sqrt{\| f\|_{\infty}}   \sqrt{c_0^2(h)c_2(h) }, 
\]

and 

\begin{eqnarray*}
\frac{Dx}{n^2} & \leq &  \frac{4  \sqrt{\| f\|_{\infty}}   \sqrt{c_0^2(h)c_2(h)}  }{ n} x\\
& \leq  & \frac{\theta' c_0^2(h) c_2(h) }{3n} + \frac{12\| f\|_{\infty}x^2}{\theta' n}.
\end{eqnarray*}

In summary, we have proved 

\begin{eqnarray*}
\frac{Ax^2}{n^2}  &\leq  & \frac{8x^2 c_0(\hmin) \| K \|_{\infty} }{n^2} \\
\frac{Bx^{3/2}}{n^2} & \leq &  \frac{\theta'}{3n} c_0^2(h) c_2(h) +  \frac{6}{\theta'} c^2_0(\hmin) c_2(\hmin) \frac{x^3}{n^2} \\
\frac{C\sqrt{x}}{n^2} & \leq &  \frac{\theta'   c_0^2(h) {c_2(h)}}{3n} +\frac{ \square  \| f\|_{\infty} x}{\theta' n} \\
\frac{Dx}{n^2} & \leq &  \frac{\theta' c_0^2(h) c_2(h) }{3n} + \frac{12\| f\|_{\infty}x^2}{\theta' n}.
\end{eqnarray*}

But 
\[ 
c_0^2(\hmin) c_2(\hmin)  = c_0(\hmin) \int  K_{\hmin^2}(x,y) c_0(\hmin)  K_{\hmin^2}(x,y) \omega_d(dy) \leq c_0(\hmin) \| K\|_{\infty}.
\] 

Thus finally with probability larger than $1-5.54 |\mathcal{H}| e^{-x}$, we have for any $h \in \mathcal{H}$ 
 \[
\left | \frac{U(h,\hmin)}{n^2}\right  |\leq   \frac{ \theta' c_0^2(h) c_2(h) }{n} + \frac{\square \| f\|_{\infty}x^2}{\theta' n} + \frac{\square c_0(\hmin)  \| K \|_{\infty}x^3} {\theta' n^2}.
\]
This ends the proof of Lemma \ref{Ustat}.

\bigskip 
Back to  (\ref{eq0}),  we have the following control 

\begin{lemma}\label{lemme11}

With probability greater that $1-9.54 | \mathcal{H} |  e^{-x}$, for any $h \in \mathcal{H}$ 
\begin{equation}
\begin{split}
\left | \langle \hat f_h -f, \hat f_{\hmin} -f \rangle -\frac{\langle c_0(h)K_{h^2}, c_0(\hmin) K_{\hmin^2} \rangle}{n} \right  | &\leq  \theta' \|f_h-f \|^2  +  \frac{ \theta' c_0^2(h) c_2(h) }{n} \\
&  \quad + \left ( \frac{\theta'}{2} + \frac{1}{2\theta'}\right) \| f_{\hmin}-f\|^2 \\
& \quad + \frac{Cx^2\| f\|_{\infty}}{\theta' n}+\frac{C(K) c_0(\hmin)x^3}{n^2},
\end{split}
\end{equation}
where $C$ is an absolute constante and $C(K)$ a constant only depending on $K$.
\end{lemma}

\textbf{Proof of Lemma \ref{lemme11}.}

We have  first to control (\ref{eq1})  and (\ref{eq2}), namely  
$$
\frac 1 n \langle  \hat f_h, f_{\hmin}\rangle -  \frac 1 n \langle  \hat f_h, f_{\hmin}\rangle  + \frac{1}{n} \langle f_h, f_{\hmin}\rangle 
$$ 
and
$$
  V(h,\hmin) +V(\hmin,h) + \langle f_h -f, f_{\hmin}-f \rangle.
 $$
 Let $h$ and $h'$ be fixed. We have
 \[
 \langle \hat f _h, f_{h'}\rangle = \frac 1 n \sum_{i=1}^n \int c_0(h)K_{h^2}(x,X_i) f_{h'}(x) \omega_d(dx).
  \]
  Therefore 
  \[
  | \langle \hat f_h, f_{h'} \rangle | \leq \| f_{h'}\|_{\infty} =  \left  \| \int c_0(h) K_{h'^2}(u,y)f(y) \omega_d(dy )  \right \|_{\infty} \leq \| f\|_{\infty},
  \]
  and 
  \[
   |\langle f_h, f_{\hmin}\rangle | \leq \| f_h\|_{\infty}  \int   f_{\hmin} (u) \omega_d(du)  \leq \| f\|_{\infty},
  \]
  which gives the control of (\ref{eq1}):
  \[
 \left |  \frac 1 n \langle  \hat f_h, f_{\hmin}\rangle -  \frac 1 n \langle  \hat f_h, f_{\hmin}\rangle  + \frac{1}{n} \langle f_h, f_{\hmin}\rangle  \right |  \leq \frac{3   \| f\|_{\infty}}{n}.
 \]
 It remains to bound the three terms of (\ref{eq2}).  We get
 
 \begin{eqnarray*}
   V(h,h') &=& \langle \hat f_h -f_h, f_{h'}-f \rangle  \\
   &=& \frac 1 n \sum_{i=1}^n \left (g_{h,h'}(X_i)-\E [g_{h,h'}(X_i)] \right )
   \end{eqnarray*}
 
 with 
 $$
 g_{h,h'}(x):= \langle c_0(h)K_{h^2}(\cdot,x), f_{h'}-f  \rangle.
 $$
 
We have
 
 $$
 \| g_{h,h'} \|_{\infty} \leq \| f_{h'}-f\|_{\infty} \leq  2 \|f \|_{\infty}.
 $$
Furthermore using Cauchy Schwarz inequality we obtain
\begin{eqnarray*}
\E [g^2_{h,h'}(X_1)] &= &\int_y  \left(  \int_x c_0(h)K_{h^2}(x,y)(f_{h'}(x)-f_{h}(x)) \omega_d(dx)\right)^2 f(y) \omega_d(dy)\\
& \leq & \| f\|_{\infty} \int_y  \left(  \int_x c_0(h)K_{h^2}(x,y)(f_{h'}(x)-f_{h}(x)) \omega_d(dx)\right)^2  \omega_d(dy) \\
&\leq &  \| f\|_{\infty} c_0(h) \int_y \left( \int_x K_{h^2}(x,y) (f_{h'}(x)-f_{h}(x))^2 \omega_d(dx) \right)  \left (c_0(h) \int_x K_{h^2}(x,y) \omega_d(dx) \right )  \omega_d(dy) \\
& \leq &  \| f\|_{\infty} c_0(h) \int_y \int_x K_{h^2}(x,y) (f_{h'}(x)-f_{h}(x))^2 \omega_d(dx)  \omega_d(dy) \\
& \leq &  \| f\|_{\infty} \| f_{h'}-f_{h}\|^2.
\end{eqnarray*}
Consequently with probability larger than $1-2e^{-x}$, Bernstein inequality  \cite{Massart} leads to 
$$
|V(h,h')| \leq \sqrt{\frac{2x}{n}  \| f\|_{\infty} \| f_{h'}-f\|^2} + \frac{2x  \|f \|_{\infty}}{3n} \\
\leq \frac{\theta'}{2}\| f_{h'}-f\|^2 + \frac{C\| f\|_{\infty}x}{\theta'n}.
$$
The bound on $V(h,h')$ obtained above  is first applied with $h'=\hmin$, then we invert the roles of $h$ and $\hmin$. Besides we have 
\[
| \langle f_h -f, f_{\hmin}-f \rangle | \leq \frac{\theta'}{2} \|f_h-f \|^2 +\frac{1}{2\theta'} \|f_{\hmin}-f \|^2.
\]

Finally using Lemma \ref{Ustat}, we get with probability larger than $1-9.54 |\mathcal{H}| e^{-x}$

\begin{eqnarray*}
\begin{split}
| \langle \hat f_h -f, \hat f_{\hmin} -f \rangle -  \frac{\langle c_0(h)K_{h^2}, c_0(\hmin) K_{\hmin^2} \rangle}{n}  | & \leq  \frac{ \theta' c_0^2(h) c_2(h) }{n} + \frac{ \square\| f\|_{\infty}x^2}{\theta' n} + \frac{\square c_0(\hmin)  \| K \|_{\infty}x^3} {\theta' n^2} \\
& \quad \quad + \frac{3   \| f\|_{\infty}}{n} +  \frac{\theta'}{2}\| f_{\hmin}-f\|^2 + \frac{\| f\|_{\infty}x}{\theta'n} + \frac{\theta'}{2}\| f_{h}-f\|^2  \\
& \quad \quad + \frac{\| f\|_{\infty}x}{\theta' n} +  \frac{\theta'}{2} \|f_h-f \|^2 +\frac{1}{2\theta'} \|f_{\hmin}-f \|^2 \\
& \leq  \theta' \|f_h-f \|^2  +  \frac{ \theta' c_0^2(h) c_2(h) }{n} +\left ( \frac{\theta'}{2} + \frac{1}{2\theta'}\right) \| f_{\hmin}-f\|^2 \\
& \quad \quad + \frac{Cx^2\| f\|_{\infty}}{\theta' n}+\frac{C(K) c_0(\hmin)x^3}{n^2},
 \end{split}
\end{eqnarray*}
which completes the proof of Lemma \ref{lemme11}.

\bigskip
Now Proposition \ref{lerasle2} gives with probability larger than $1-\square |\mathcal{H}| e^{-x}$, for any $h \in \mathcal{H}$:
\[
\| f-f_h\|^2 + \frac{ c_0^2(h)c_2(h)}{n}  \leq  2 \| f -\hat f_h\|^2 + C_2(K) \frac{\| f\|_{\infty} x^2}{ n },
\]
where $C_2(K)$ depends only on $K$. Hence by applying Lemma \ref{lemme11} with $h$ first and then $\hat h$ we obtain with probability larger than $1-\square |\mathcal{H}| e^{-x}$, for any $h \in \mathcal{H}$:
\begin{equation}\label{eq5}
\begin{split}
& \left   |  \langle \hat f_h -f, \hat f_{\hmin} -f  \rangle -\frac{\langle c_0(h)K_{h^2},  c_0(\hmin)K_{\hmin^2} \rangle}{n}  -  \langle \hat f_{\hat h} -f, \hat f_{\hmin} -f  \rangle +\frac{\langle c_0(\hat h) K_{\hat h^2},  c_0(\hmin) K_{\hmin^2} \rangle }{n}  \right | \\
& \leq  2\theta' \| \hat f_h -f\|^2 +2\theta' \| \hat f_{\hat h } -f\|^2 + \left ( \theta'+\frac{1}{\theta'} \right ) \|f_{\hmin}-f \|^2 +\frac{\tilde C(K)}{\theta'} \left(\frac{\| f\|_{\infty} x^2}{n} +\frac{x^3c_0(\hmin)}{n^2 } \right),
\end{split}
\end{equation}
where $\tilde C(K)$ is a constant only depending on $K$. 
Now back to (\ref{16}) and using (\ref{eq5}), we have
\begin{equation*}
\begin{split}
\| \hatfhat -f\|^2 & \leq \| \hat f_h -f\|^2  \\
 & \quad  +   pen_\lambda(h)  -2 \left ( \langle \hat f_h -f,  \hat f_{\hmin}-f \rangle -   \frac{\langle c_0(h)K_{h^2},  c_0(\hmin)K_{\hmin^2} \rangle}{n} \right  ) -  2\frac{\langle c_0(h)K_{h^2},  c_0(\hmin)K_{\hmin^2} \rangle}{n} \\
&\quad -    pen_\lambda(\hat h) + 2 \left(  \langle \hatfhat -f, \hat f_{\hmin}- f \rangle   -   \frac{\langle c_0(h)K_{\hat h^2},  c_0(\hmin)K_{\hmin^2} \rangle}{n}  \right )  +  2\frac{\langle c_0(\hat h)K_{\hat h^2},  c_0(\hmin)K_{\hmin^2} \rangle}{n} \\
& \leq  \| \hat f_h -f\|^2  + pen_\lambda(h) -2 \frac{\langle c_0(h)K_{h^2},  c_0(\hmin)K_{\hmin^2} \rangle}{n}  -    pen_\lambda(\hat h)  + 2 \frac{\langle c_0(\hat h)K_{\hat h^2},  c_0(\hmin)K_{\hmin^2} \rangle}{n} \\
& \quad  + 4   \theta' \| \hat f_h -f\|^2 + 4\theta' \| \hat f_{\hat h } -f\|^2 + 2 \left ( \theta'+\frac{1}{\theta'} \right ) \|f_{\hmin}-f \|^2 +\frac{\tilde C(K)}{\theta'} \left(\frac{\| f\|_{\infty} x^2}{n} +\frac{x^3c_0(\hmin)}{n^2 } \right),
\end{split}
\end{equation*}
 choosing $\theta'=\frac \theta 4 $ yields the result. This completes the proof of Proposition \ref{theorem9}.

\bigskip

\textbf{Proof of Theorem \ref{theorem 2}.}

We set $\tau= \lambda-1$. Let $ \varepsilon \in(0,1)$ and $\theta \in (0,1)$ depending on $\varepsilon$ to be specified later.  
Developing the expression of  $\pen_\lambda(h)$ given in (\ref{penalite}), we  have that 
\begin{eqnarray*}
\pen_\lambda(h) &= & \frac{\lambda c_0^2(h) c_2(h)}{n}-\frac{c_0^2(\hmin)c_2(\hmin)}{n} -  \frac{c_0^2(h)c_2(h)}{n} +2 \frac{\langle c_0(h) K_{h^2} , c_0(\hmin)K_{\hmin^2} \rangle }{n} \\
&=& \tau \frac{c_0^2(h) c_2(h)}{n} -\frac{c_0^2(\hmin)c_2(\hmin)}{n} + 2 \frac{\langle c_0(h) K_{h^2} , c_0(\hmin)K_{\hmin^2} \rangle }{n}.
\end{eqnarray*}

Using Proposition \ref{theorem9} and the expression of  $\pen_\lambda(h)$ given above, we obtain with probability greater than $1-\square |\mathcal{H}| \exp(-x)$, for any $h \in \mathcal{H}$,

\begin{equation}\label{preuvetheo2eq1}
\begin{split}
(1-\theta)\| \hat f_{\hat h} -f\|^2 +\tau \frac{c_0^2(\hat h)c_2(\hat h)}{n} & \leq (1+\theta) \| \hat f_h -f\|^2 + \tau \frac{c_0^2(h)c_2(h)}{n} +\frac{C_2}{\theta}  \| f_{\hmin} -f\|^2  \\
&\quad \quad + \frac{C(K)}{\theta} \left (  \frac{\| f\|_{\infty}x^2}{n} +\frac{c_0(\hmin)x^3}{n^2}\right).
\end{split}
\end{equation}
We first consider the case $\tau \geq 0$.  Using  (\ref{lerasle22}) of Proposition  \ref{lerasle2},  with probability $1- \square |\mathcal{H}| e^{-x}$ one has
\[
\tau \frac{c_0^2(h)c_2(h)}{n} \leq \tau(1+\theta) \| f-\hat f_h\|^2 + \tau  \frac{C'(K) \|f \|_{\infty} x^2}{\theta^3 n}, 
\]
where $C'(K)$ is a constant only depending on the kernel $K$. As $\tau \frac{c_0^2(\hat h)c_2(\hat h)}{n} \geq 0 $, thus (\ref{preuvetheo2eq1}) becomes \begin{equation}
\begin{split}
(1-\theta) \| \hat f_{\hat h} -f\|^2  & \leq (1+\theta+\tau(1+\theta)) \| \hat f_h -f\|^2 +\frac{C_2}{\theta}\| f_{\hmin} -f\|^2  \\
&\quad + \frac{C(K)}{\theta} \left ( \frac{\| f\|_{\infty}x^2}{n} + \frac{c_0(\hmin)x^3 }{n^2}\right) +\tau \frac{C'(K)\| f\|_{\infty}x^2}{\theta^3n }.
\end{split}
\end{equation}

With $\theta = \varepsilon/(\varepsilon + 2 + 2\tau)$, we obtain
\begin{equation*}
\begin{split}
\| \hat f_{\hat h} -f\|^2  \leq (1+\tau + \varepsilon) \| \hat f_{ h} -f\|^2  &+\frac{C_2(\varepsilon +2 + 2\tau)^2}{(2+2\tau)\varepsilon}  \| f_{\hmin} -f\|^2 \\
&+C''(K, \varepsilon, \tau) \left( \frac{\| f\|_{\infty}x^2}{n} + \frac{c_0(\hmin)x^3 }{n^2}\right),
\end{split}
\end{equation*}
with $C''(K, \varepsilon, \tau)$ a constant depending only on $K, \varepsilon, \tau$.

Let us now study the case $-1 < \tau \leq 0$. Using  (\ref{lerasle22}) of Proposition \ref{lerasle2} with $h=\hat h$, we have with probability  $1-\square |\mathcal{H}| e^{-x}$
$$
\tau \frac{c_0^2(\hat h)c_2(\hat h)}{n} \geq \tau(1+\theta) \| f-\hat f_{\hat h}\|^2 + \tau \frac{C'(K)\| f\|_{\infty}x^2}{\theta^3 n},
$$
with $C'(K)$ a constant depending only on $K$.  Consequently, as $\tau \frac{c_0^2(h)c_2(h)}{n} \leq 0$, (\ref{preuvetheo2eq1}) becomes
\begin{equation*}
\begin{split}
(1-\theta + \tau(1+\theta))\| \hat f_{\hat h} -f \|^2 & \leq (1+\theta) \| \hat f_{h} -f \|^2 +\frac{C_2}{\theta}  \| f_{\hmin} -f\|^2  \\
& + \frac{C(K)}{\theta} \left ( \frac{\| f\|_{\infty}x^2}{n}  +  \frac{c_0(\hmin) x^3}{n^2}\ \right) - \tau \frac{C'(K)\| f\|_{\infty}x^2}{\theta^3 n}.
\end{split}
\end{equation*}
With $\theta=(\varepsilon(\tau+1)^2)/(2+\varepsilon(1-\tau^2))<1$, we obtain with probability $1- \square |\mathcal{H}| e^{-x}$,
\begin{equation*}
\begin{split}
\| \hat f_{\hat h} -f \|^2  \leq \left( \frac{1}{1+\tau}  + \varepsilon \right) \| \hat f_{ h} -f \|^2 &+ C''(\varepsilon, \tau)   \| f_{\hmin} -f\|^2  \\
& + C'''(K, \varepsilon, \tau) \left( \frac{\| f\|_{\infty}x^2}{n} + \frac{c_0(h_{\min})x^3}{n^2 }\right).
\end{split}
\end{equation*}
This completes the proof of Theorem \ref{theorem 2}.

\bigskip

\textbf{Proof of Theorem \ref{theorem_petit}}.

We still set $\tau= \lambda -1$.  We set $\theta \in (0,1)$ such that $\theta < -(1+\tau)/5$. We consider  inequality (\ref{preuvetheo2eq1}) written with $h=\hmin$. One obtains 

\begin{equation*}
\begin{split}
(1-\theta)\| \hat f_{\hat h} -f\|^2 +\tau \frac{c_0^2(\hat h)c_2(\hat h)}{n} & \leq (1+\theta) \| \hat f_{\hmin} -f\|^2 + \tau \frac{c_0^2(\hmin)c_2(\hmin)}{n} +\frac{C_2}{\theta}  \| f_{\hmin} -f\|^2  \\
&\quad \quad + \frac{C(K)}{\theta} \left (  \frac{\| f\|_{\infty}x^2}{n} +\frac{c_0(\hmin)x^3}{n^2}\right).
\end{split}
\end{equation*}

Now consider equation (\ref{lerasle21}) with $h=\hmin$, one gets
\[
\| f -\hat f_{\hmin} \|  \leq (1+\theta)  \left ( \|f-f_{\hmin} \|^2  +\frac{c_0^2(\hmin)c_2(\hmin) }{n}\right) +  \frac{C'(K) \| f\|_{\infty}· x^2}{\theta^3 n}.
\]
Combining the two inequalities above, we have
\begin{equation}
\begin{split}
(1-\theta)\| \hat f_{\hat h} -f\|^2 +\tau \frac{c_0^2(\hat h)c_2(\hat h)}{n} &\leq  \left ((1+\theta)^2 + \frac{C_2}{\theta} \right) \|f-f_{\hmin} \|^2 \\
& +  (\tau +(1+\theta)^2) \frac{c_0^2(\hmin)c_2(\hmin) }{n} \\
&+ \frac{C(K)}{\theta} \left  (   \frac{ \| f\|_{\infty}· x^2}{ n} + \frac{ c_0(\hmin) x^3}{n^2}  \right )  \\
& + (1+\theta) \frac{C'(K) \|f \|_{\infty}x^2}{\theta^3 n }.
\end{split}
\end{equation}

Now let us define $u_n:=\frac{ \| f_{\hmin} -f\|^2}{c_0^2(\hmin) c_2(\hmin)/n}$. We have by assumption that $u_n \rightarrow 0$ when $n$ tends to infinity. Then we get,
\begin{equation}\label{5.18}
\begin{split}
(1-\theta) \| \hat f_{\hat h} -f\|^2 +\tau \frac{c_0^2(\hat h)c_2(\hat h)}{n}  & \leq \left  ( \left (1+\theta)^2 + \frac{C_2}{\theta}  \right )u_n + \tau + (1+ \theta)^2 \right ) \frac{c_0^2(\hmin) c_2(\hmin)}{n} \\
& \quad +  C(K, \theta) \left  (  \frac{ \| f\|_{\infty}· x^2}{ n}  + \frac{x^3}{n^2}c_0(\hmin) \right ).
\end{split}
\end{equation}

No we consider equation (\ref{lerasle22}) with $h= \hat  h$ and $\eta=1$, we get
\[
\frac{c_0(\hat h) c_2(\hat h)}{n} \leq 2 \| f -\hat f_{\hat h}\| + C'(K) \frac{\| f\|_{\infty}x^2}{n},
\]
then
\[
\| f -\hat f_{\hat h}\| \geq \frac{c_0(\hat h) c_2(\hat h)}{2n} - C'(K) \frac{\| f\|_{\infty}x^2}{n}. 
\]

And (\ref{5.18}) becomes 
\[
\begin{split}
\left (\frac{1-\theta}{2}+ \tau \right ) \frac{c_0^2(\hat h) c_2(\hat h)}{n} & \leq \left  ( \left ( (1+ \theta)^2 + \frac{C_2}{\theta} \right ) u_n + \tau + (1+ \theta)^2 \right ) \frac{c_0^2(\hmin)c_2(\hmin) }{n} \\
& \quad  + C'(K,\theta) \left ( \frac{\| f\|_{\infty} x^2}{n}  + \frac{x^3 c_0(\hmin)}{n^2}\right ).
\end{split}
\]
But we assumed that $u_n=o(1)$. Thus for $n$ large enough $((1+ \theta)^2+ \frac {C_2}{\theta})u_n \leq \theta.$
We are now going to bound the remaining terms $  C'(K,\theta)\left ( \frac{\| f\|_{\infty} x^2}{n}  + \frac{ c_0(\hmin) x^3}{n^2}\right )$. We have
\[
\begin{split}
 C'(K,\theta)\left ( \frac{\| f\|_{\infty} x^2}{n}  + \frac{ c_0(\hmin)x^3}{n^2}\right ) \frac{n}{c_0^2(\hmin)c_2(\hmin)} &= C''(K, \theta, \|f\|_{\infty})\bigg  (\frac{x^2}{c_0^2(\hmin)c_2(\hmin)} \\
  &\quad \quad  + \frac{x^3}{n c_0(\hmin) c_2(\hmin)} \bigg ) \\
 &\leq  C''(K, \theta, \|f\|_{\infty}) \left(  x^2 \hmin ^{d-1}+ \frac{x^3}{n} \right ),
\end{split}
\]
for $n$ large enough using  ($\ref{c01}$) and (\ref{c21}). 
But $\hmin^{d-1} \leq \frac{(\log n )^\beta }{n}$ and setting $x=( \frac{n}{\log n})^{\frac 1 3}$, we get
\begin{eqnarray*}
C''(K, \theta, \|f\|_{\infty}) \left(  x^2 \hmin ^{d-1}+ \frac{x^3}{n} \right ) \leq C''(K, \theta, \|f\|_{\infty}) \left ( \frac{(\log n)^{\beta -2/3}}{n^{1/3}} + \frac{1}{\log n} \right ) =o(1 )\leq \theta, 
\end{eqnarray*}
for $n$ large enough. Consequently there exists $N$ such that for $n \geq N$, with probability larger than $1-\square |\mathcal{H}| e^{-(n/\log n)^{1/3}}$
\[
\left (\frac{1-\theta}{2}+ \tau \right ) \frac{c_0^2(\hat h) c_2(\hat h)}{n} \leq (\theta + \tau +(1+\theta)^2 + \theta)\frac{ c_0^2(\hmin)c_2(\hmin)}{n} \leq (1+ \tau +5\theta) \frac{ c_0^2(\hmin)c_2(\hmin)}{n}.
\]
 Using ($\ref{c0c2}$) we have for $n$ large enough that
\[
0.9 \,  h^{1-d}R(K) \leq c_0^2(h)c_2(h) \leq 1.1 \,h^{1-d} R(K)
\]

and thus we finally get for $n$ large enough 
\[
0.9 \left (\frac{1-\theta}{2}+ \tau \right ) \hat h^{1-d} \leq 1.1 (1+ \tau +5\theta) \hmin^{1-d}   \iff   \left (\frac{1-\theta}{2}+ \tau \right ) \hat h^{1-d} \leq 1.23(1+ \tau +5\theta) \hmin^{1-d}.
\]
But $\frac{(1-\theta)}{2} + \tau < 1+ \tau < 0$, and because we have chosen $\theta $ such that  $ 1+ \tau +5\theta <0$ (for instance $\theta = -(\tau+1)/10$)), one gets
\[
\hat h \leq \left(  \frac{1.23(1+ \tau +5\theta) }{\frac{1-\theta}{2}+ \tau} \right )^{\frac{1}{d-1}}\hmin.
\]
With  $\theta = -(\tau+1)/10$, the inequality above becomes for $n$ large enough
\[
\hat h \leq  \left( 1.23 \left( 2.1-\frac 1 \lambda \right ) \right)^{\frac{1}{d-1}}\hmin,
\]
which completes the proof of Theorem \ref{theorem_petit}. 

\bigskip 

\textbf{Proof of Theorem \ref{vitesse}}.

As usual, the MISE of $\hat f_{\check h }$ is divided into a bias term and a variance term.  Let us focus on the bias term. Define for $f : \S \rightarrow \R$ and $s$ even,
$$
\mathcal{D}^s f= \sum_{i=1}^{s/2} \frac{2i}{(2i)!} \gamma_{2i,s/2-i} D^{2i}f,
$$ 
where 
$$   \gamma_{i} = \sum_{\alpha_1+\cdots + \alpha_{d-1}=i} \frac{(-1)^{\alpha_1}}{(2\alpha_1+1)!} \cdots  \frac{(-1)^{\alpha_{d-1}}}{(2\alpha_{d-1}+1)!}, \quad \gamma_0=1. $$
A bound for the bias is given in  \cite{Klemela2}. We  stated it in the next Proposition.
\begin{proposition}\label{convergenceconvolution} Assume that $f \in \textbf{F}_2(s)$. Let $K$ be a class $s$ kernel, where $s \geq 2$ is even. Then
$$\lim_{h \rightarrow 0 } \| h^{-s} |  \E(\hat f_h) -f | - | \alpha_0^{-1}(K) \alpha_s(K) \mathcal{D}^s f | \| =0.$$
\end{proposition}

Now let $f \in  \tilde{\textbf{F}}_2(s, B)$ and $\mathcal{E}$ the event corresponding to the intersection of events in Corollary \ref{corollaire} and Proposition \ref{lerasle2}. For any $A>0$, by taking $x$ proportional to $\log n $, $ \P(\mathcal{E} )\geq 1-n^{-A}$. We have
\begin{eqnarray*}
\E \left [\| \hat f_{\check h} -f \|^2 \right ] \leq \E \left [\| \hat f_{\check h} -f \|^2 \mathds{1}_{\mathcal{E} }  \right  ] +  \E \left [\| \hat f_{\check h} -f \|^2 \mathds{1}_{\mathcal{E}^c }  \right  ].
\end{eqnarray*}
Let us deal with the second term of the right hand side.  We have
$$
  \| \hat f_{ h} -f \|^2   \leq 2( \| \hat f_h\|^2 + \| f\|^2).
$$
But
\begin{eqnarray*}
 \| \hat f_h\|^2& =& \frac{c^2_0(h)}{n^2} \sum_{i,j} \int_{\S}  K_{h^2}(x,X_i) K_{h^2}(x,X_j) \omega_d(dx) \\
 & \leq & \frac{c_0(h)}{n^2} \| K\|_{\infty} \sum_{i,j } c_0(h) \int_{\S}  K_{h^2}(x,X_j) \omega_d(dx) \\
 & \leq &  c_0(h) \| K\|_{\infty} \leq 2n,
\end{eqnarray*}
since  $c_0(h) \leq 2 n/ \|  K\|_{\infty}$,  using (\ref{c01}) and (\ref{H}) for $n$ large enough.  Thus 
$$
  \| \hat f_{ h} -f \|^2 \leq 2n  + 2 \| f\|^2,
$$
which gives the result on  $ \mathcal{E}^c$.
Now on $ \mathcal{E}$, for $n\geq n_0$ ($n_0$ not depending on $f$) Proposition \ref{lerasle2} and  Proposition \ref{convergenceconvolution}  yield that 
\begin{eqnarray*}
\min_{h \in \mathcal{H}}  \| \hat f_{ h} -f \|^2  &\leq & (1+ \eta) \min_{h \in \mathcal{H}}  \left( \|f-f_h \|^2  +\frac{c_0^2(h)c_2(h) }{n} \right )+ \square \frac{\Upsilon (\log n) ^2}{\eta^3 n}  \\
& \leq & C   \min_{h \in \mathcal{H}}  \left (h^{2s} + \frac{h^{1-d}}{n} \right ) + \square \frac{\Upsilon (\log n) ^2}{\eta^3 n}.
\end{eqnarray*}
Minimizing in $h$ the right hand side of the last inequality gives the result on $ \mathcal{E}$.

\section{Appendix}\label{appendix}

Next lemma collects some standard properties about constants $c_0$ and $c_2$ that are useful when dealing with kernel density estimation with directional data. 
\begin{lemma}\label{proprietes}
We have 
\begin{equation}
\label{c01}
c_0^{-1}(h)= R_0(K)h^{d-1} +o(1),
\end{equation}
as $h \rightarrow 0$ uniformly in $h$, we recall that $ R_0(K)= 2^{(d-3)/2}\sigma_{d-1}\alpha_0(K)$.
We also have 
\begin{equation}
\label{c21}
c_2(h)=R_1(K)h^{d-1}+o(1),
\end{equation}
as $h \rightarrow 0$, uniformly in $h$, with $R_1(K)$ defined in (\ref{R1}). Eventually we have
\begin{equation}
\label{c0c2}
c_0^2(h)c_2(h)=h^{1-d}R(K)+o(1),
\end{equation}
when $h \rightarrow 0$   uniformly in $h$, with
$$R(K):= R_1(K)/R_0^2(K).$$
\end{lemma}
The proof of Lemma \ref{proprietes} can be found in the proof of Proposition 4.1 of  \cite{Amiri}. 

\bigskip
In practice, SPCO algorithm turns to be simple to compute as  shown in the next lemma.  
\begin{lemma}\label{SPCO-pratique}
For $\mathbb{S}^2$ and $K(x)=e^{-x}$, we have  that 
\begin{eqnarray*}
\| \hat f_h -\hat f_{\hmin}\|^2 &=&   \frac{4\pi c^2_0(h)}{n^2} e^{-2/h^2} h^2 \sum_{i,j} \frac{\sinh ( |X_i+X_j| /h^2)}{|X_i+X_j|}  \\
&\quad  + &  \frac{4\pi c^2_0(h_{\min})}{n^2} e^{-2/{h^2_{\min}}} h^2_{\min} \sum_{i,j} \frac{\sinh ( |X_i+X_j| /h^2_{\min})}{|X_i+X_j|} \\
& \quad  - & \frac{8 \pi }{n^2} c^2_0(h) c^2_0(h_{\min}) e^{-1/{h^2}} e^{-1/{h^2_{\min}}}  \sum_{i,j}\frac{\sinh ( |X_i/h^2 +X_j/h^2_{\min}| )}{ |X_i/h^2 +X_j/h^2_{\min}| },
\end{eqnarray*}
with
\begin{eqnarray*}
c_0(h)^{-1}=4\pi e^{-1/h^2} h^2 \sinh(1/h^2),
\end{eqnarray*}
and
\[
c_2(h) = 2\pi e^{-2/h^2}h^2\sinh(2/h^2).
\]
\end{lemma}
\textbf{Proof of Lemma \ref{SPCO-pratique}}.

We have 
\begin{eqnarray*}
\| \hat f_h -\hat f_{\hmin}\|^2 &=& \| \hat f_h \|^2 + \| \hat f_{\hmin} \|^2 -2 \int \hat f_h \hat f_{\hmin},
\end{eqnarray*}
but
\begin{eqnarray*}
 \| \hat f_h \|^2 & =&  \int_{\mathbb{S}^2} \left |  \frac{c_0(h)}{n} \sum_i e^{-(1-x^TX_i)/h^2} \right |^2 \omega_d(dx) \\
 &=&  \frac{c^2_0(h)}{n^2} \sum_{i,j}  \int_{\mathbb{S}^2} e^{-(1-x^TX_i)/h^2} e^{-(1-x^TX_j)/h^2} \omega_d(dx) \\
 &=&  \frac{c^2_0(h)}{n^2} e^{-2/h^2} \sum_{i,j}  \int_{\mathbb{S}^2} e^{ x^T \frac{X_i+X_j}{|X_i+X_j|} \frac{ |X_i+X_j|}{h^2}} \omega_d(dx) \\
 &=& \frac{2\pi c^2_0(h)}{n^2} e^{-2/h^2} \sum_{i,j} \int_{-1}^{1}e^{\frac{|X_i+X_j|}{h^2}t}dt \\
 &=&  \frac{4\pi c^2_0(h)}{n^2} e^{-2/h^2} h^2 \sum_{i,j} \frac{\sinh ( \frac{|X_i+X_j|} {h^2})}{|X_i+X_j|}.
\end{eqnarray*}
Hence following the same computations, one finds that 
\begin{eqnarray*}
\| \hat f_h -\hat f_{\hmin}\|^2 &=&   \frac{4\pi c^2_0(h)}{n^2} e^{-2/h^2} h^2 \sum_{i,j} \frac{\sinh ( |X_i+X_j| /h^2)}{|X_i+X_j|}  \\
&\quad  + &  \frac{4\pi c^2_0(h_{\min})}{n^2} e^{-2/{h^2_{\min}}} h^2_{\min} \sum_{i,j} \frac{\sinh ( |X_i+X_j| /h^2_{\min})}{|X_i+X_j|} \\
& \quad  - & \frac{8 \pi }{n^2} c^2_0(h) c^2_0(h_{\min}) e^{-1/{h^2}} e^{-1/{h^2_{\min}}}   \sum_{i,j}\frac{\sinh ( |X_i/h^2 +X_j/h^2_{\min}| )}{ |X_i/h^2 +X_j/h^2_{\min}| }.
\end{eqnarray*}
Furthermore, using Remark \ref{remarque1}
\begin{eqnarray*}
c_0(h)^{-1}=   \int_{\mathbb{S}^2} e^{-(1-x^Ty)/h^2} dy = 2\pi \int_{-1}^1 e^{-(1-t)/h^2}dt =4\pi e^{-1/h^2} h^2 \sinh(1/h^2),
\end{eqnarray*}
and
\[
c_2(h) =  \int_{\mathbb{S}^2} e^{-2(1-x^Ty)/h^2} dy= 2\pi \int_{-1}^1  e^{-2(1-t)/h^2}dt= 2\pi e^{-2/h^2}h^2\sinh(2/h^2).
\]
This ends the proof of Lemma \ref{SPCO-pratique}. 

\vspace{1.5cm}

\textbf{Acknowledgments} The author would like to thank Claire Lacour and Vincent Rivoirard for interesting suggestions and remarks.

\end{document}